\documentclass[a4paper,12pt]{article} 
\usepackage{amsmath}                                  
\usepackage{amssymb}                                  
\usepackage{graphicx}                                 
\usepackage{booktabs}                                 
\usepackage[english]{babel}                           
\usepackage[a4paper,                                  
            top = 25mm,                               
            bottom = 25mm,                            
            left = 25mm,                              
            right = 25mm]{geometry} 
\usepackage[dvipsnames]{xcolor}
\usepackage{fancyhdr}
\usepackage{listings}
\usepackage{amsthm}
\usepackage{mathtools}
\usepackage{subfigure}
\usepackage{hyperref}
\usepackage{marvosym}
\usepackage{chemformula}
\newcommand\numberthis{\addtocounter{equation}{1}\tag{\theequation}}
\usepackage[utf8]{inputenc}
\usepackage{wrapfig}   
\usepackage{multicol}
\begin{document}                                    
\fancyfoot{}
\fancyfoot[C]{\thepage}

\thispagestyle{plain}
\begin{center}
    \large
    \textbf{Optimal Control of a Power Storage Facility with Variable Payoffs}
        
    \vspace{0.4cm}
    \textbf{Fraser J W O'Brien\footnote{Maxwell Institute for Mathematical Sciences}\textsuperscript{,}\footnote{Department of Actuarial Mathematics and Statistics, Heriot–Watt University, Edinburgh EH14 4AS, UK}, $\;$ Timothy C Johnson\footnotemark[\value{footnote}]}

    \normalsize
    
\end{center}

\begin{abstract}
    We present a methodology for determining the relationship between the optimal control points of a power storage facility and a number of different factors including storage level and temperature. The interaction between different factors is considered to allow for the identification of a precise optimal control strategy to maximise the profits of a power storage facility under a variety of different conditions.  The methodology is based upon traditional stochastic techniques, however by working directly with excess demand data is model independent and does not require identification of the underlying process.
\end{abstract}

\section{Introduction}

In order to combat the effects of climate change, the UK is transitioning away from its reliance on polluting fossil fuels to a more renewable energy focused grid. Whilst this transition comes with many environmental benefits, it also introduces logistical problems in dealing with the increased variability in power generation due to renewable energy sources' dependence on natural conditions. The higher levels of variability leads to increased difficulty in balancing the supply and demand of electricity, potentially resulting in the price of electricity going negative as supply outstrips demand. Whilst initially this may seem beneficial, negative prices occur at times of peak productivity for renewable energy plants, reducing income and causing concern for the overall profitability of renewable energy schemes. Such concerns could discourage investment in more sustainable sources of energy, hampering attempts for carbon-neutrality.

Power storage facilities are able to convert excess electricity on the grid into other forms of energy for storage until needed, buying electricity cheaply when an excess is produced and selling for profit as demand increases. With large enough total storage capacity, storage facilities can prevent the price of electricity going negative by ensuring that demand is always able to meet supply. Power storage facilities support the profitability of renewable schemes, whilst bringing down costs for the general public during times of low renewable generation. The EU has estimated that total energy storage capacity will need to rise more than three-fold from 2022 to 2030 to accommodate the increasing supply of renewable energy \cite{EU}. Pumped storage hydro accounts for almost all of the power storage capacity within the UK, however with limited geographically suitable locations remaining alternative forms of power storage will need to be considered to provide flexibility. An overview of the properties and potential uses of different forms of power storage can be found in \cite{B-I}.

\newpage

Any private investment in power storage technologies is made with a consideration of the profitability. As such, the theory within this paper considers the perspective of potential investors rather than the grid. Grid operators, policy makers and regulators must consider different sets of criteria when discussing the optimal control of power storage facilities, as they aim to provide a low-cost and flexible power grid which is capable of meeting the demands of the nation. However through the use of models such as those seen within this paper, policy makers can gain an understanding of the effect potential market mechanisms may have on facility operating procedures.

The price of electricity in the UK is determined by the consideration of the merit-order curve (MOC), an approximation of which can be found in Figure 1. Power plants submit bids with the quantity of electricity they can provide and for what cost, which are then sorted based on price to determine the collection of plants required to meet the range of demand values with the cheapest being prioritised. One effect of this system is that the price of electricity can exhibit jumps as costlier forms of electricity generation become necessary. Whilst the price of electricity depends on the demand, this relationship depends on the amount of electricity being generated through renewable sources. Non-renewable sources of electricity are controllable and thus their total capacity is available at all times (excluding instances where power plants are down for maintenance), ensuring a relatively consistent relationship between non-renewable sources of electricity and prices. Power plants (or storage facilities) are paid relative to the bid submitted by the most expensive power plant included to meet the demand from the nation, with cheaper forms of electricity generation rewarded for their cost effectiveness with increased profit margins.

\begin{figure}[h!]
    \centering
    \includegraphics[width=0.9\textwidth]{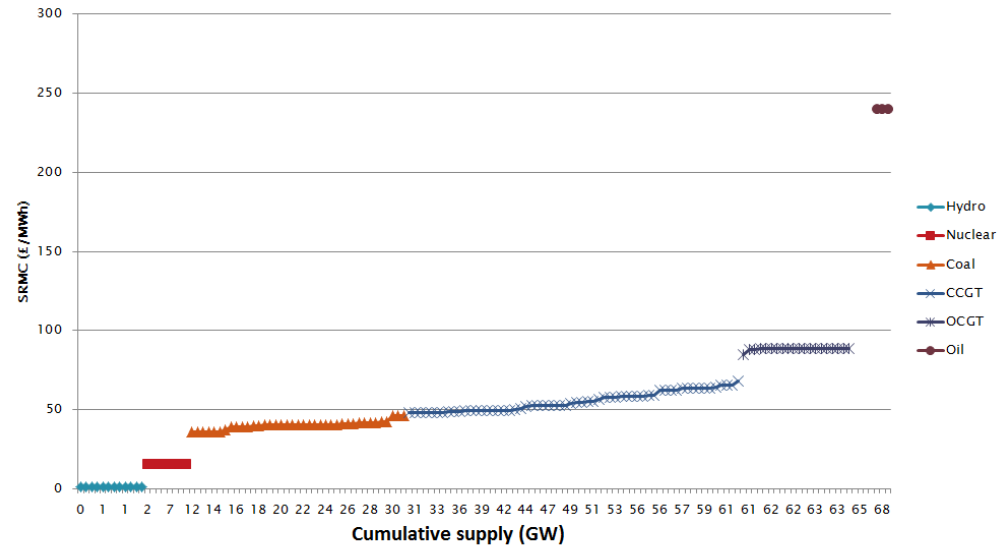}
    \caption{The MOC approximation for the UK market excluding renewable energy sources \cite{EG}.}
\end{figure}

\newpage

Note that the inclusion of renewable energy would result in a translation in the cumulative supply, as renewables have low running costs. At this point we introduce the concept of excess demand as the difference between the demand for electricity and supply of renewable electricity, and due to the variance in power output of renewable energy sources and the consistency of other forms can be considered a more appropriate variable for the price of electricity than demand. Due to the stochasticity in the output of renewable energy sources, excess demand can be assumed to be generated by an underlying stochastic process. Control of power storage facilities can be considered either an optimal stopping or optimal switching problem depending on the type of facility considered. A review of the potential approaches to optimise energy management can be found in \cite{M-M}.

Of interest to this paper is an optimal switching problem, however similarities in the theory enables the use of results from papers focusing on impulse control \cite{A},\cite{C} and optimal stopping \cite{D-K}, including \cite{J} which significantly relaxes a number of assumptions in order to provide a more realistic and flexible environment to consider the optimal stopping problem. Different approaches to tackle these problems include Markovian \cite{A} and Martingale approaches \cite{D-Z}, however both of these methodologies require the identification of the underlying process generating the excess demand data for applicable results. A novel approach was introduced in \cite{EG} \cite{EG-J}, in which the optimal control strategy is estimated directly from excess demand data through the estimation of the values of waiting without the need to estimate the stochastic process assumed to underlie the data. Such an approach is model independent and as such can suitably be applied to data where the form of the underlying process is unknown.

It is worthwhile for facility operators to employ a model which considers the problem of maximising the profit for the specific facilities they operate, and to consider how the physical properties of the facility may effect the optimal control strategy. \cite{C-PD} consider the optimal control strategy of variable speed pumped storage hydropower plants consisting of Francis turbines that participate in secondary regulation reserve markets, indicating the level of precision and number of considerations necessary to identify the optimal control strategy for an individual facility. The methodology presented here is a development of the work of \cite{EG} \cite{EG-J}, in the latter of which the authors admit "the identified strategy is naive, involving a single stopping region for each regime", and the aim of this paper is to expand the sophistication of the optimal control strategy by considering the effect of a number of different factors, such that optimal control schemes can be tailored to individual facilities from the same excess demand data-sets. Whilst one option to account for certain factors is to consider a multi-dimensional diffusion generating excess demand (such processes have been considered in works such as \cite{C-I}), this paper also aims to provide a methodology for which sophisticated mathematical knowledge is not required.

\newpage
\section{Theoretical Basis}
The facility operator wishes to charge and discharge the facility in such a way to maximise the profits made by the facility. This can be expressed in mathematical terms by the desire to maximise the following performance criterion:

$$J(\mathbb{Z}_{z,x}) := \lim_{n \rightarrow \infty} \mathbb{E}_x \left[ \sum_{j=1}^{n-1} e^{-\Lambda_{T_j}}\left[ E(X_            {T_j})\boldsymbol{1}_{\{\Delta Z_{T_j} < 0\}} - F(X_{T_j})\boldsymbol{1}_{\{\Delta Z_{T_j} > 0\}} \right]\boldsymbol{1}_{\{T_j < \infty\}} \right]$$

where $\mathbb{Z}_{z,x}$ represents a chosen control strategy, $X_t$ the level of excess demand, at time $t$, $E$ the income gained by discharging the facility, $F$ the cost of charging the facility, $T_j$ selected switching times, $Z_t$ the storage level of the facility, and $\Lambda_t = \int_{0}^{t} r(X_s) \,ds$ is a state-dependant factor by which all future pay-offs are discounted. For the scenarios considered within this paper, we take $r(x) = r$ for some constant $r > 0$. The aim is to identify the control strategy which provides the value function $v$ that maximises this performance criterion, given by:
$$v(x, z) = \sup_{\mathbb{Z}_{z,x} \in \mathcal{A}_{z,x}} J(\mathbb{Z}_{z,x}) \quad \text{ for } x \in \mathcal{I}, \quad z \in Z $$
Where $\mathcal{A}_{z,x}$ is the set of all admissible strategies. To ensure that this optimisation problem is well posed, we assume:
$$E(x)-F(x) < 0 \quad \text{ for all } x \in \mathcal{I}$$
to prevent the instantaneous charging and discharging of the facility for profit. This is a suitable assumption to make as the round-trip efficiency of all possible power storage facilities is less than 100\%. We can consider the value of our different options through the following expressions, if we act and claim the payoff then the value of our choice is given by:
\begin{align*}
    v(x, z) \ge
          \begin{dcases}
            - F(x) \quad &\text{ if } \Delta z > 0 \\
            E(x) \quad &\text{ if } \Delta z < 0  
            \end{dcases}
\end{align*}

Alternatively the facility can wait a short amount of time $\Delta t$ before proceeding optimally, giving an expected reward of:
$$v(x, z) \ge \mathbb{E}^x[e^{-\Lambda_{\Delta t}}v(X_{\Delta t}, z)]$$

The option which provides the most value for the facility is the optimal choice at any time, meaning that we need to compare the immediate value of acting against the potential reward of waiting a short time before proceeding optimally. By waiting we hope that the price of electricity will become more favourable, however this must be balanced against the potential loss of value of the facility due to the discount factor. Since the immediate payoffs are known quantities, we must estimate the expected value of waiting in order to determine the optimal strategy.

\newpage

The expected value of waiting depends on the properties of the excess demand data, which we assume to be a stochastic process driven by some one-dimensional It\^{o} diffusion given by the following:

\begin{equation}
    dX_t = b(X_t)dt + \sigma(X_t)dW_t, \quad X_0 = x_0 \in \mathcal{I} 
\end{equation}

where $\mathcal{I}$ is an open interval with end points $-\infty < \alpha \le X_t \le \beta < \infty$ and $W_t$ is standard Brownian motion. By taking the expected value for waiting, applying It\^{o}'s formula, integrating between $0$ and $\Delta t$ then dividing by $\Delta t$ and taking the limit $\Delta t \rightarrow 0$ provides the equation:

\begin{equation}
    \frac{1}{2}\sigma^2(x)\frac{\partial^2 v}{\partial x^2}(x, z) + b(x)\frac{\partial v}{\partial x}(x, z) - rv(x, z) \le 0
\end{equation}

the solution to the equality gives the expected value of waiting before proceeding optimally. Depending upon the form of the process underlying the excess demand data (and thus the form of $\sigma^2(x)$ and $b(x)$) this equation can be solved analytically. However the underlying process for the excess demand data is not known, or may even not have a canonical form. As such, we must find an alternative approach to identify the solution to the differential equation. Following the work of \cite{J-Z}, we begin by making the following assumptions:

Assumption 1: The functions $b$, $\sigma: \mathcal{I} \rightarrow \mathbb{R}$ are $\mathcal{B}(\mathcal{I})$-measurable with:
$$\sigma^2(x) > 0, \text{ for all } x \in \mathcal{I}$$
and 
$$\int_{a}^{b} \frac{1 + |b(s)|}{\sigma^2(s)} ds < \infty \text{ and } \sup_{s \in [\alpha, \beta]} \sigma^2(s) < \infty, \text{ for all }  \alpha < a < b < \beta.$$

Assumption 2: The solution to (1) is non-explosive.

Note that Assumption 3 in \cite{J-Z} is satisfied through the choice of constant $r > 0$, and that the underlying process does not need to be known for these assumptions to hold, as is the case with excess demand data. Using these assumptions ensures that the general solution to (2) exists and is given by 
$$v(x, z) = A(z)\phi(x) + B(z)\psi(x)$$
for constants with respect to $x$, $A(z), B(z) \in \mathbb{R}$, with the functions $\phi, \psi$ defined between two levels of excess demand $x, y$ as follows:
$$
    \mathbb{E}_x(e^{-r\tau_{y}}) = 
          \begin{dcases}
            \frac{\psi(x)}{\psi(y)} \quad \text{ if } x \le y\\
            \frac{\phi(x)}{\phi(y)} \quad \text{ if } x \ge y
            \end{dcases}
$$

\newpage
which can be rearranged to the following:
\begin{align*}
\psi(y) = 
    \begin{dcases}
        \psi(x) \mathbb{E}_y \left[ e^{-r\tau_{x}} \right] \quad &\text{ if } x \le y \\
        \frac{\psi(x)}{\mathbb{E}_y\left[ e^{-r\tau_{x}} \right]} \quad &\text{ if } x \ge y
    \end{dcases}
\end{align*}
\begin{align*}
\phi(y) = 
    \begin{dcases}
        \frac{\phi(x)}{\mathbb{E}_y\left[ e^{-r\tau_{x}} \right]} \quad &\text{ if } x \le y \\
        \phi(x) \mathbb{E}_y \left[ e^{-r\tau_{x}} \right] \quad &\text{ if } x \ge y 
    \end{dcases}
\end{align*}
These functions have absolutely continuous first derivatives, and are unique up to the multiplicative constants. Furthermore as $\tau_x > 0$, we can see the following properties:
\begin{align*}
    0 < \psi(x) \quad &\text{and} \quad \psi'(x) > 0, \quad \text{ for all } x \in I \\
    0 < \phi(x) \quad &\text{and} \quad \phi'(x) < 0, \quad \text{ for all } x \in I \\
    \text{and} \quad \lim_{x \uparrow \beta} \psi(x) = &\lim_{x \downarrow \alpha} \phi(x) = \infty.
\end{align*}

These definitions provide a direct link between the excess demand data and the value of waiting, as the general solution to (2) can be estimated by calculating the expected hitting times between levels of  excess demand $x, y$ across the whole range of excess demand values. To ensure that the switching problem remains well-posed, we require the following assumption:

Assumption 3: The payoff functions $E, F: \mathcal{I} \rightarrow \mathbb{R}$ are the difference of two convex functions, and we have
$$\lim_{x \uparrow \beta} \frac{|F(x)|}{\psi(x)} = \lim_{x \downarrow \alpha} \frac{|F(x)|}{\phi(x)} = 0$$
$$\lim_{x \uparrow \beta} \frac{|E(x)|}{\psi(x)} = \lim_{x \downarrow \alpha} \frac{|E(x)|}{\phi(x)} = 0$$
whilst the limits $\lim_{x \downarrow \alpha} \frac{F(x)}{\phi(x)}, \lim_{x \uparrow \beta} \frac{F(x)}{\psi(x)}, \lim_{x \downarrow \alpha} \frac{E(x)}{\psi(x)}$ and $ \lim_{x \uparrow \beta} \frac{E(x)}{\phi(x)}$ exist in $[-\infty, \infty]$. This assumption guarantees the existence of a region in which it is optimal to act.

The properties of $\phi(x), \psi(x)$ are well-stated and can be found in \cite{B-S}, \cite{F}, \cite{I-H} for example.

\newpage
\section{Optimal Control Points of a Multi-Use facility}

Power storage facilities capable of sequentially charging and discharging are of particular interest in this paper, as they are necessary to provide flexibility to a renewable energy focus grid. There are however alternative forms of power storage requiring different models which may play a part in facilitating the construction of renewable energy sources, such as in \cite{L-EA} in which a battery is located onsite at a wind farm, charging when electricity is not being sold to the grid. We first review the work of \cite{EG} and consider the so called "bang-bang" case, in which a facility can only be in states $z \in [0,1]$ representing either fully discharged or fully charged respectively. In such a scenario, the value functions for the facility are given by:
\begin{align*}
    w(x,0) = 
          \begin{dcases}
            A_1\phi(x) + B_1\psi(x) - F(x), &\text{ if } x \in (\alpha, a] \\
            A_0\phi(x) + B_0\psi(x), &\text{ if } x \in (a, \beta),
            \end{dcases}
\end{align*}
\begin{align*}
    w(x,1) = 
          \begin{dcases}
            A_1\phi(x) + B_1\psi(x), &\text{ if } x \in (\alpha, b] \\
            A_0\phi(x) + B_0\psi(x) + E(x), &\text{ if } x \in (b, \beta),
            \end{dcases}
\end{align*}

At this point it is useful to note that the value of $B(z)\phi(x)$ can be thought of as representing the value of the charge of the facility, and $A(z)\psi(x)$ the value of the remaining empty capacity of the facility. Thus we can deduce that $A_1 = B_0 = 0$ giving the following:
\begin{align*}
    w(x,0) = 
          \begin{dcases}
            B\psi(x) - F(x), &\text{ if } x \in (\alpha, a] \\
            A\phi(x), &\text{ if } x \in (a, \beta),
            \end{dcases}
\end{align*}
\begin{align*}
    w(x,1) = 
          \begin{dcases}
            B\psi(x), &\text{ if } x \in (\alpha, b] \\
            A\phi(x) + E(x), &\text{ if } x \in (b, \beta),
            \end{dcases}
\end{align*}

where the subscript has been removed for ease of notation. The boundaries between the continuation and switching regions $(a,b)$ can be identified by equating sides of the value function and using the principle of smooth fit, as seen in \cite{EG}. The result of such a derivation is that we aim to identify the points $(a,b)$ that satisfy:
\begin{align}
    A &= \frac{d}{dx} \left( \frac{F(x)}{\psi(x)} \right)_{x=a} \frac{\psi^2(a)}{\mathcal{W}(a)} = \frac{d}{dx} \left( \frac{E(x)}{\psi(x)} \right)_{x=b} \frac{\psi^2(b)}{\mathcal{W}(b)} \\
    B &= \frac{d}{dx} \left( \frac{F(x)}{\phi(x)} \right)_{x=a} \frac{\phi^2(a)}{\mathcal{W}(a)} = \frac{d}{dx} \left( \frac{E(x)}{\phi(x)} \right)_{x=b} \frac{\phi^2(b)}{\mathcal{W}(b)}.
\end{align}

Where $\mathcal{W}(x) = \phi(x)\psi'(x) - \phi'(x)\psi(x)$ denotes the Wronskian. By identifying the functions $\phi$ and $\psi$, it is possible to calculate the optimal control points $(a,b)$ without the need to identify the underlying process. Algorithm 2 in \cite{EG} allows for identification of the optimal control points from this system of equations.

\newpage
\section{Marginal Control of Power Storage Facilities}
We now consider a facility which is capable of partial charges and discharges, allowing for more control over the buying and selling of electricity. By considering a facility acting in a marginal sense rather than a 'bang-bang' scenario allows for a more realistic model for the control of a power storage facility, and allows for a more detailed control strategy to be identified. As such there are additional considerations to be made, with the optimal control points potentially varying with storage level.

Where previously the storage facility could only obtain two states, we now have $z \in [z_0, z_1, ... , z_n]$, where $z_0 = 0\%$ capacity, $z_n = 100\%$, and for simplicity we assume that $z_1 - z_0 = z_2 - z_1 \text{ } ... \text{ } z_n - z_{n-1}$. We can represent the changing efficiency of our storage facility by including $z$ dependence in our payoff functions, with $E(x), F(x)$ becoming $E(x,z), F(x,z)$, now representing the payoff of switching to the next level of storage rather than full/empty. We re-introduce the subscript notation for our constants $A(z_i) = A_i, B(z_i) = B_i$, and the optimal control points $a(z_i) = a_i, b(z_i) = b_i$. The aim now becomes to maximise the following performance criterion:
$$J(\mathbb{Z}_{z,x}) := \lim_{n \rightarrow \infty} \mathbb{E}_x \left[ \sum_{j=1}^{n-1} e^{-\Lambda_{T_j}}\left[ E(x_{T_j}, z_{T_j})\boldsymbol{1}_{\{ \Delta z_{T_j} < 0\}} - F(x_{T_j}, z_{T_j})\boldsymbol{1}_{\{ \Delta z_{T_j} > 0\}} \right]\boldsymbol{1}_{\{T_j < \infty\}} \right].$$

Beginning with an empty storage facility, we are presented with the choice to either charge the facility or wait before proceeding optimally. As such the value function is similar to the previous case, given by:
\begin{align*}
    w(x,z_0) = 
          \begin{dcases}
            A_1\phi(x) + B_1\psi(x) - F(x,z_0), &\text{ if } x \in (\alpha, a_0] \\
            A_0\phi(x), &\text{ if } x \in (a_0, \beta),
            \end{dcases}
\end{align*}

And similarly with a full facility, we have:
\begin{align*}
    w(x,z_n) = 
          \begin{dcases}
            B_n\psi(x), &\text{ if } x \in (\alpha, b_n] \\
            A_{n-1}\phi(x) + B_{n-1}\psi(x) + E(x,z_n), &\text{ if } x \in (b_n, \beta).
            \end{dcases}
\end{align*}

Recalling that $A_n = B_0 = 0$. Of particular interest for this scenario is deriving the optimal control strategy for storage values $z_i, 0 < i < n$, which can be found by considering the following value function:
\begin{align*}
    w(x,z_i) = 
          \begin{dcases}
            A_{i+1}\phi(x) + B_{i+1}\psi(x) - F(x,z_i), &\text{ if } x \in (\alpha, a_i] \\
            A_i\phi(x) + B_i\psi(x), &\text{ if } x \in (a_i, b_i) \\
            A_{i-1}\phi(x) + B_{i-1}\psi(x) + E(x,z_i), &\text{ if } x \in [b_i, \beta) 
            \end{dcases}
\end{align*}

To ease notation, we now define $F_i(x) = F(x,z_i), E_i(x)=E(x,z_i)$. For a full derivation of the following result, see Appendix A. It can be found that the above value functions result in the identification of the pair of optimal control points $(a_i, b_{i+1})$, which satisfy the following system of equations:

\begin{align}
    -\left[ A_{i+1}-A_{i} \right] &= \frac{d}{dx}\left( \frac{F_i(x)}{\psi(x)} \right)_{x = a_i} \frac{\psi^2(a_i)}{\mathcal{W}(a_i)} = \frac{d}{dx}\left( \frac{E_{i+1}(x)}{\psi(x)} \right)_{x = b_{i+1}} \frac{\psi^2(b_{i+1})}{\mathcal{W}(b_{i+1})} \\
    \left[ B_{i+1}-B_{i} \right] &= \frac{d}{dx}\left( \frac{F_i(x)}{\phi(x)} \right)_{x = a_i} \frac{\phi^2(a_i)}{\mathcal{W}(a_i)} = \frac{d}{dx}\left( \frac{E_{i+1}(x)}{\phi(x)} \right)_{x=b_{i+1}} \frac{\phi^2(b_{i+1})}{\mathcal{W}(b_{i+1})}
\end{align}

for $0 \le i \le n-1$. Thus it is possible to identify the optimal control points for the full operating range of the power storage facility.

The system of equations (5-6) can be solved in a similar manner to (3-4) in order to identify the pair of optimal control points $(a_i, b_{i+1})$. Note that $a_i$ is the point where it becomes optimal to charge from storage level $z_i$ to $z_{i+1}$, whilst $b_{i+1}$ is the point where it becomes optimal to discharge from storage level $z_{i+1}$ to $z_i$. The coefficients $A_i, B_i$ can be found by identifying the values of $A_{n-1}, B_1$ and the differences $\left[ A_{i}-A_{i+1} \right], \left[ B_{i+1}-B_{i} \right], 0 \le i \le {n-1}$. 

\section{Variation of the Optimal Control Points}

There are two types of factors to consider when calculating the optimal control strategy of a power storage facility, those which affect the demand for electricity and those which affect the storage facility of interest. When working with realistic payoff functions, the price of electricity varies depending upon certain thresholds crossed by excess demand, thus when considering factors affecting the whole grid rather than just the single power storage facility the factors can not be separated from the level of excess demand. However factors affecting the amount of electricity required to act in a certain manner (typically related to the physical properties of the storage facility) can be separated from the price of electricity, as the total cost of acting is given by the price of electricity multiplied by the amount of electricity required. 

\subsection{Variation with Storage Level}

Of interest is the effect of the variation of the payoff functions $F(x, z), E(x, z)$ with respect to $z$ on the optimal control scheme. One possible method to identify the full optimal control strategy is to solve equations (5-6) for the full range of values of $z$, however this may prove computationally expensive depending on the number of values. 

\newpage

At this point we establish the necessary functions for use in simulations throughout the remainder of this paper. The results obtained within this paper consider a known Ornstein-Uhlenbeck process given by the following: 
\begin{align*}
    dX_t &= \kappa (\theta - X_t)dt + \sigma dW_t \\
    \kappa &= 0.003\\
    \theta &= 5,000 \\
    \sigma &= 900. 
\end{align*}

By defining the process underlying the excess demand data, the problems of working with real data are avoided, whilst still providing a demonstration of the potential applications of the algorithms in \cite{EG}. Whilst Geometric Brownian motion has been used in works such as \cite{D-Z}, \cite{S} found that the use of mean-reverting processes is preferable and Ornstein-Uhlenbeck processes have seen use in \cite{C-L} and \cite{EG}. The above values for the Ornstein-Uhlenbeck process reflect a future scenario with a higher percentage of renewable energy, in which excess demand can range from $-40,000$MW during times of low demand and high renewable, to $50,000$MW at times of high demand and minimal renewable generation. The fundamental solutions $\phi, \psi$ for an Ornstein-Uhlenbeck process as defined above are given by: 
$$\psi(x) = e^{\frac{\kappa (x-\theta)^2}{2\sigma^2}} D_{-\frac{r}{\kappa}}\left( \sqrt{\frac{2\kappa}{\sigma}}(x-\theta) \right)$$
$$\phi(x) = e^{\frac{\kappa (x-\theta)^2}{2\sigma^2}} D_{-\frac{r}{\kappa}}\left( \sqrt{\frac{2\kappa}{\sigma}}(\theta-x) \right)$$
where $D_v$ is the parabolic cylinder function with index $v$. 

We consider linear payoff functions with respect to excess demand which result in the largest effects on the optimal control points by additional factors, given by the following:
\begin{align*}
    F_x(x) &= 0.001x + 20 \\
    E_x(x) &= 0.9 \times (0.001x + 20)
\end{align*}

where the parameter values are chosen to approximate the gradient of the MOC for Coal and CCGT, and allow for negative price values when a surplus of renewable energy is being generated. Work by \cite{J-Z} and \cite{L-Z} allows for pay-off functions not in $C^2$, allowing for realistic payoff functions to be defined including price jumps, however for the sake of simplicity the complexity introduced by price jumps is omitted. The multiplication by $0.9$ can be considered as representing a facility with $90\%$ round-trip efficiency.

As defined previously, $Z$ denotes the storage state of the storage facility, and we define functions $Z_f(z), Z_e(z)$ which represent the effect of the storage level on the payoff functions $F, E$ respectively. Depending on the type of facility, the overall inefficiencies of the system are increased as the energy in the system increases, and as such we define the following general functions:
\begin{align*}
    Z_f(z) &= 1 + 2\frac{z}{z_n}\\
    Z_e(z) &= 1
\end{align*}
which correspond to a change in efficiency from $90\%$ to $30\%$, recalling $z_n$ is the maximum capacity of the storage facility. Whilst individual facilities may not see variation this large, this wider range will include those covered by individual facilities, as well showcasing the range of optimal control strategies used by facilities depending upon their physical properties. Therefore the payoff functions can be defined by:
\begin{align*}
    F(x, z) &= F_x(x)Z_f(z) = (0.001x + 20)\left( 1 + 2\frac{z}{z_n} \right) \\
    E(x, z) &= E_x(x)Z_e(z) = 0.9 \times (0.001x + 20)
\end{align*}
With these definitions, the relationship between the optimal control points $(a,b)$ and the storage level $z$ can be seen in Figure 2:

\begin{figure}[h!]
    \centering
    \includegraphics[width=0.75\textwidth]{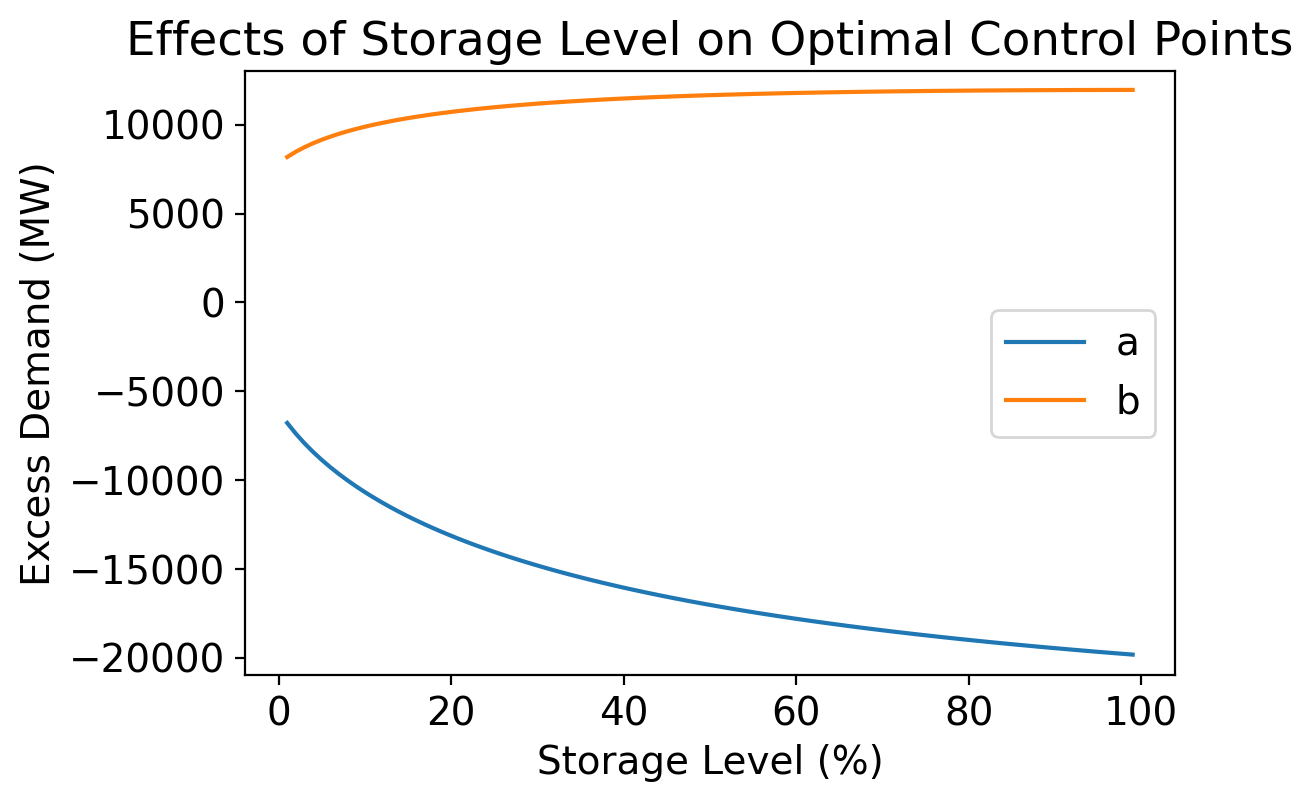}
    \caption{The effect of storage level on the optimal control points.}
\end{figure}

It can be seen that as storage level increases (and therefore the inefficiency in the system), the optimal control points move into more extreme regions of excess demand. In particular, the lower optimal control point $a$ decreases from $-6787.10$MW to $-19826.88$MW. This is to be expected, as the increased inefficiencies in the system lead to higher buying prices for electricity, meaning the value of waiting before acting increases pushing the boundary $a$ into more extreme values of excess demand. This increased value of waiting also affects $b$ to a lesser extent, increasing from $8183.25$MW to $11961.30$MW, due to the increased expected time between acting allowing for additional waiting for the selling price to increase.

Whilst it is possible to calculate the optimal control points explicitly at every level of storage $z$, it would be preferable to quantify the relationship between the optimal control points and $z$. Note that whilst we consider $z$ in the following analysis, the results hold for all separable variables $y$ (and can be seen by a similar approach replacing $Z_f(z_{i}), Z_e(z_{i+1}), Z_f(z_{i+1}), Z_e(z_{i+2})$ with $Y_f(y_{i}), Y_e(y_{i}), Y_f(y_{i+1}), Y_e(y_{i+1})$ respectively, and with the pairing $(a_i, b_{i+1})$ replaced by $(a_i, b_i)$ to denote the level of the variable $y$ rather than the switching between levels of storage).

\newpage
We proceed by considering the level of storage $z$, assuming that the payoff functions are separable and can be defined by:
\begin{align*}
    F(x, z) &= F_x(x)Z_f(z) \\
    E(x, z) &= E_x(x)Z_e(z)
\end{align*}

Thus the system of equations (5-6) becomes:
\begin{align}
    \left[ A_{i}-A_{i+1} \right] &= Z_f(z_i)\frac{d}{dx}\left( \frac{F_x(x)}{\psi(x)} \right)_{x=a_i} \frac{\psi^2(a_i)}{\mathcal{W}(a_i)} = Z_e(z_{i+1})\frac{d}{dx}\left( \frac{E_{x}(x)}{\psi(x)} \right)_{x=b_{i+1}} \frac{\psi^2(b_{i+1})}{\mathcal{W}(b_{i+1})} \\
    \left[ B_{i+1}-B_{i} \right] &= Z_f(z_i)\frac{d}{dx}\left( \frac{F_x(x)}{\phi(x)} \right)_{x=a_i} \frac{\phi^2(a_i)}{\mathcal{W}(a_i)} = Z_e(z_{i+1})\frac{d}{dx}\left( \frac{E_{x}(x)}{\phi(x)} \right)_{x=b_{i+1}} \frac{\phi^2(b_{i+1})}{\mathcal{W}(b_{i+1})}
\end{align}

For ease of notation we now make the following definitions:
\begin{align*}
    q_{\psi, f}(x) &= \frac{d}{dx}\left( \frac{F_x(x)}{\psi(x)} \right) \frac{\psi^2(x)}{\mathcal{W}(x)}\\
    q_{\psi, e}(x) &= \frac{d}{dx}\left( \frac{E_x(x)}{\psi(x)} \right) \frac{\psi^2(x)}{\mathcal{W}(x)}\\
    q_{\phi, f}(x) &= \frac{d}{dx}\left( \frac{F_x(x)}{\phi(x)} \right) \frac{\phi^2(x)}{\mathcal{W}(x)}\\
    q_{\phi, e}(x) &= \frac{d}{dx}\left( \frac{E_x(x)}{\phi(x)} \right) \frac{\phi^2(x)}{\mathcal{W}(x)}
\end{align*}

such that we can re-write the system of equations (7-8) as:
\begin{align}
    \left[ A_{i}-A_{i+1} \right] &= Z_f(z_i)q_{\psi, f}(a_i) = Z_e(z_{i+1})q_{\psi, e}(\ell^A(a_i)) \\
    \left[ B_{i+1}-B_{i} \right] &= Z_f(z_i)q_{\phi, f}(a_i) = Z_e(z_{i+1})q_{\phi, e}(\ell^B(a_i))
\end{align}

By considering the values of $(a_{i+1}, b_{i+2})$, found by solving the following system of equations:
\begin{align}
    \left[ A_{i+1}-A_{i+2} \right] &= Z_f(z_{i+1})q_{\psi, f}(a_{i+1}) = Z_e(z_{i+2})q_{\psi, e}(\ell^A(a_{i+1})) \\
    \left[ B_{i+2}-B_{i+1} \right] &= Z_f(z_{i+1})q_{\phi, f}(a_{i+1}) = Z_e(z_{i+2})q_{\phi, e}(\ell^B(a_{i+1}))
\end{align}

given knowledge of the optimal control points $(a_i, b_{i+1})$, it can be found that the derivatives of the optimal control points with respect to the variable $z$ are given by the following equations:
\begin{align*}
    \frac{da}{dz} &= \left[ \frac{Z'_e(z)}{Z_e(z)} - \frac{Z'_f(z)}{Z_f(z)} \right] \left( \frac{q_{\psi,f}(a)q'_{\phi,e}(b) - q'_{\psi,e}(b)q_{\phi,f}(a)}{q'_{\psi,f}(a)q'_{\phi,e}(b) - q'_{\psi,e}(b)q'_{\phi,f}(a)} \right) \\
    &= \frac{d}{dz}\left( \frac{Z_e(z)}{Z_f(z)} \right)\left( \frac{q_{\psi,e}(b)q'_{\phi,e}(b) - q'_{\psi,e}(b)q_{\phi,e}(b)}{q'_{\psi,f}(a)q'_{\phi,e}(b) - q'_{\psi,e}(b)q'_{\phi,f}(a)} \right) \numberthis
\end{align*}
\begin{align*}
    \frac{db}{dz} &= \left[ \frac{Z'_e(z)}{Z_e(z)} - \frac{Z'_f(z)}{Z_f(z)} \right] \left( \frac{q_{\psi,e}(b)q'_{\phi,f}(a) - q'_{\psi,f}(a)q_{\phi,e}(b)}{q'_{\psi,f}(a)q'_{\phi,e}(b) - q'_{\psi,e}(b)q'_{\phi,f}(a)} \right) \\
    &= \frac{d}{dz}\left( \frac{Z_f(z)}{Z_e(z)} \right)\left( \frac{q_{\psi,f}(a)q'_{\phi,f}(a) - q'_{\psi,f}(a)q_{\phi,f}(a)}{q'_{\psi,f}(a)q'_{\phi,e}(b) - q'_{\psi,e}(b)q'_{\phi,f}(a)} \right) \numberthis
\end{align*}

For the sake of brevity, a complete derivation of these results is omitted here, but given in Appendix B. Unfortunately it is not possible to integrate these derivatives to find a continuous relationship between the optimal control points $(a, b)$ and the variable $z$. This is due to the fact that there is no explicit relationship between the two optimal control points $a, b$, or the derivatives of $q_{\psi,f}(a), q_{\psi, e}(b)$ and $q_{\phi,f}(a), q_{\phi, e}(b)$. As such it is impossible to integrate the above expressions in such a way that does not cancel back down to the initial equations (9-10). However this result is still useful, as it allows for numerical methods to be used to calculate the optimal control strategy across the range of $z$ values, increasing the computational speed for a more detailed scheme to be identified than through explicit calculation for all $z$ values from the equations (9-10).

In the derivations of equations (13-14) a finite difference approximation for the gradients of $a,b$ can be found, given by the following:

\small
\begin{align*}
    \frac{\Delta a(z_i)}{\Delta z} = \frac{1}{Z_f(z_i) + \Delta z Z'_f(z_i)}\Bigg( &Z'_e(z_{i+1})\left[ \Bigg(\frac{ q_{\psi,e}(b(z_{i+1}))q'_{\phi,e}(b(z_{i+1})) - q'_{\psi,e}(b(z_{i+1}))q_{\phi,e}(b(z_{i+1}))}{q'_{\psi,f}(a(z_i))q'_{\phi,e}(b(z_{i+1})) - q'_{\psi,e}(b(z_{i+1}))q'_{\phi,f}(a(z_i))} \Bigg) \right] \\
    + &Z'_f(z_i)\Bigg[  \Bigg(\frac{ q'_{\psi,e}(b(z_{i+1}))q_{\phi,f}(a(z_i)) - q_{\psi,f}(a(z_i))q'_{\phi,e}(b(z_{i+1})) }{q'_{\psi,f}(a(z_i))q'_{\phi,e}(b(z_{i+1})) - q'_{\psi,e}(b(z_{i+1}))q'_{\phi,f}(a(z_i))} \Bigg) \Bigg] \Bigg) \numberthis
\end{align*}

\begin{align*}
        \frac{\Delta b(z_{i+1})}{\Delta z} = \frac{1}{Z_e(z_{i+1}) + \Delta z Z'_e(z_{i+1})}\Bigg( &Z'_e(z_{i+1})\left[\frac{ q_{\psi,e}(b(z_{i+1}))q'_{\phi,f}(a(z_i)) - q'_{\psi,f}(a(z_i))q_{\phi,e}(b(z_{i+1})) }{q'_{\psi,f}(a(z_i))q'_{\phi,e}(b(z_{i+1})) - q'_{\psi,e}(b(z_{i+1}))q'_{\phi,f}(a(z_i))}\right] \\
        + &Z'_f(z_i)\left[\frac{ q'_{\psi,f}(a(z_i))q_{\phi,f}(a(z_i)) - q_{\psi,f}(a(z_i))q'_{\phi,f}(a(z_i)) }{q'_{\psi,f}(a(z_i))q'_{\phi,e}(b(z_{i+1})) - q'_{\psi,e}(b(z_{i+1}))q'_{\phi,f}(a(z_i))}\right] \Bigg). \numberthis
\end{align*}

\normalsize
to find the optimal control points $(a_{i+1},b_{i+2})$ given the points $(a_i,b_{i+1})$ at. We can find the optimal control points $(a_0, b_1)$ through the use of the algorithm described in \cite{EG}. The results of this numerical approach can be seen in Figure 3.

\begin{figure}[h!]
	\centering
	\subfigure{\includegraphics[width=0.45\textwidth]{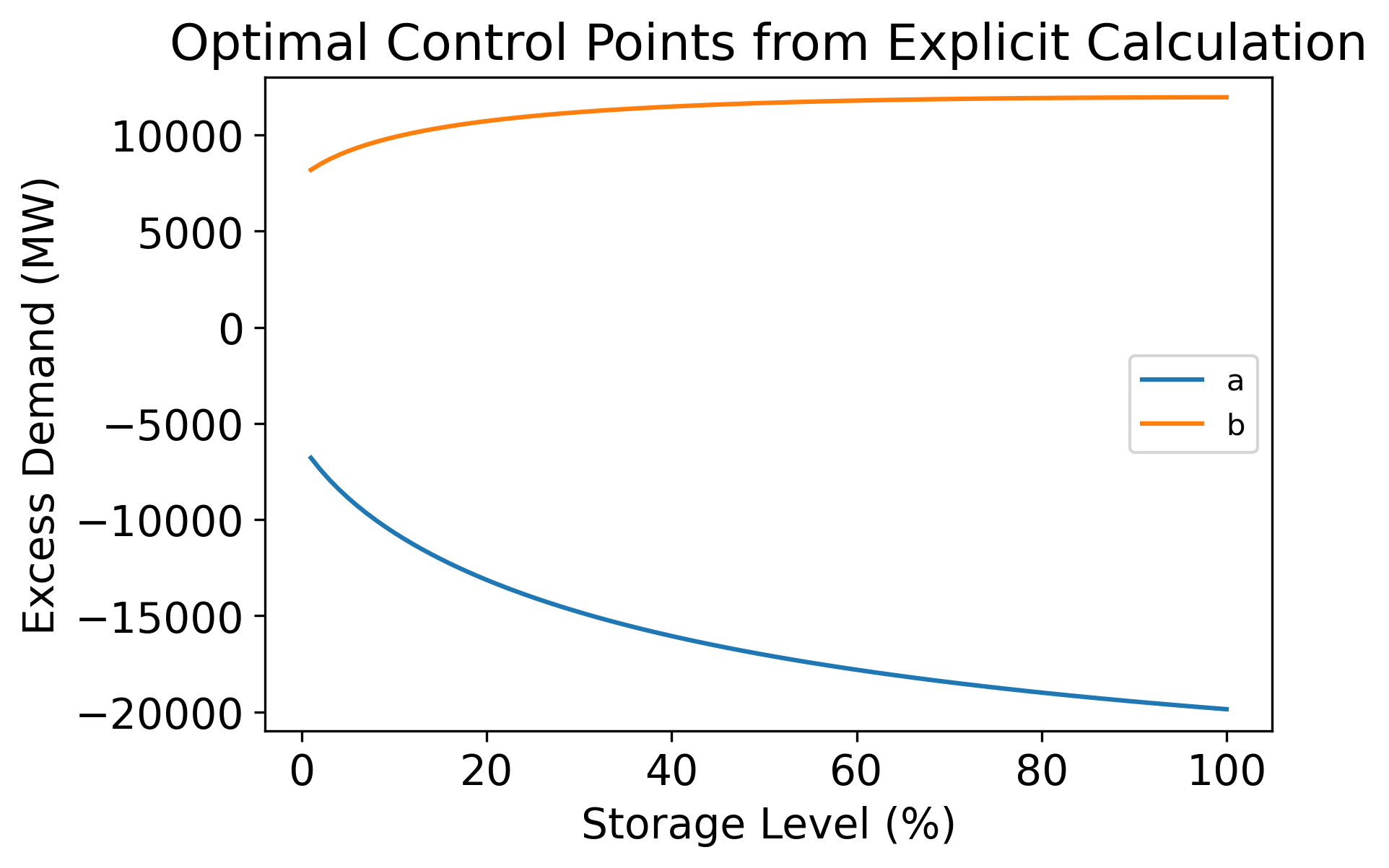} } 
	\subfigure{\includegraphics[width=0.45\textwidth]{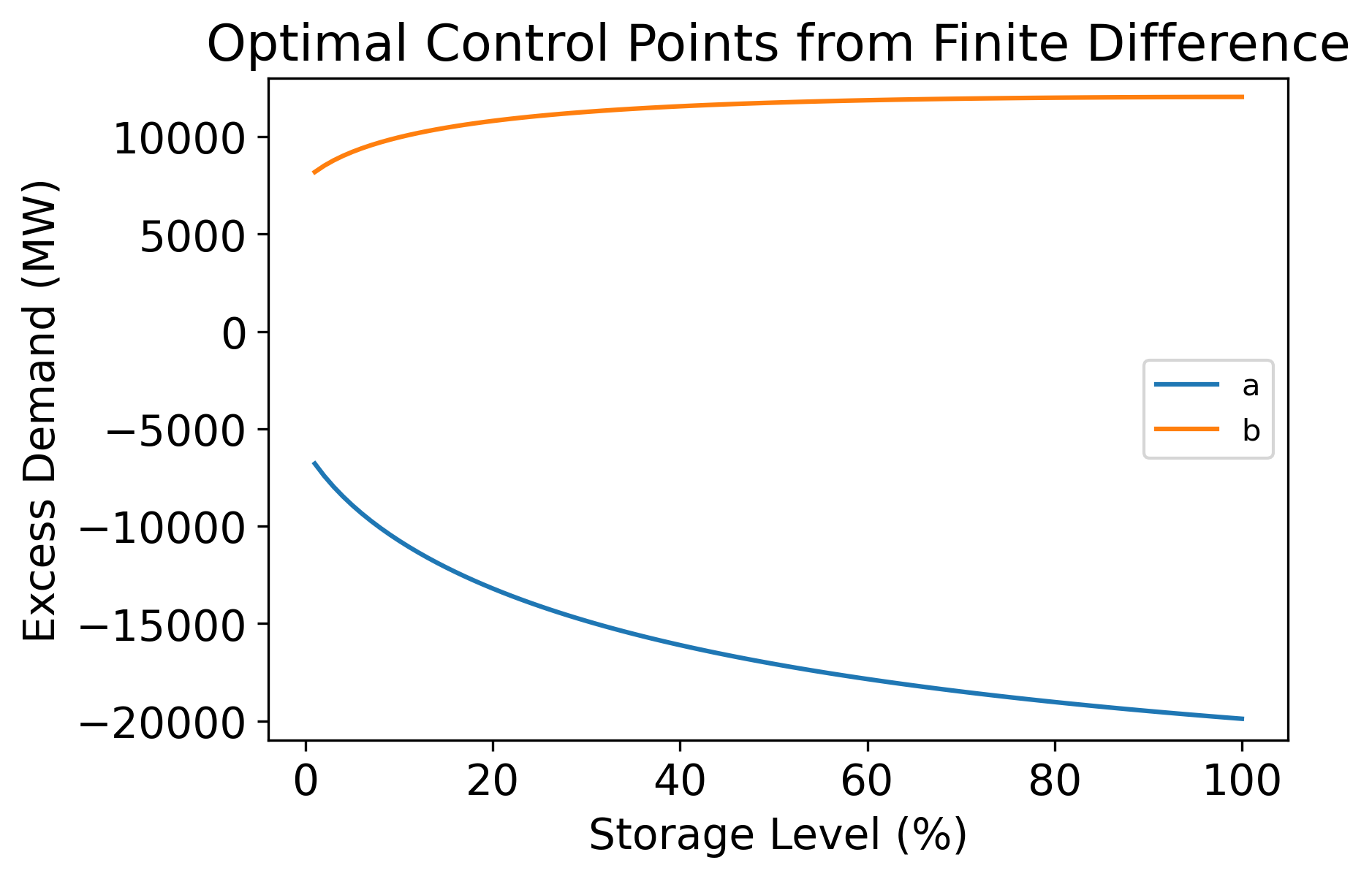} }
    \subfigure{\includegraphics[width=0.45\textwidth]{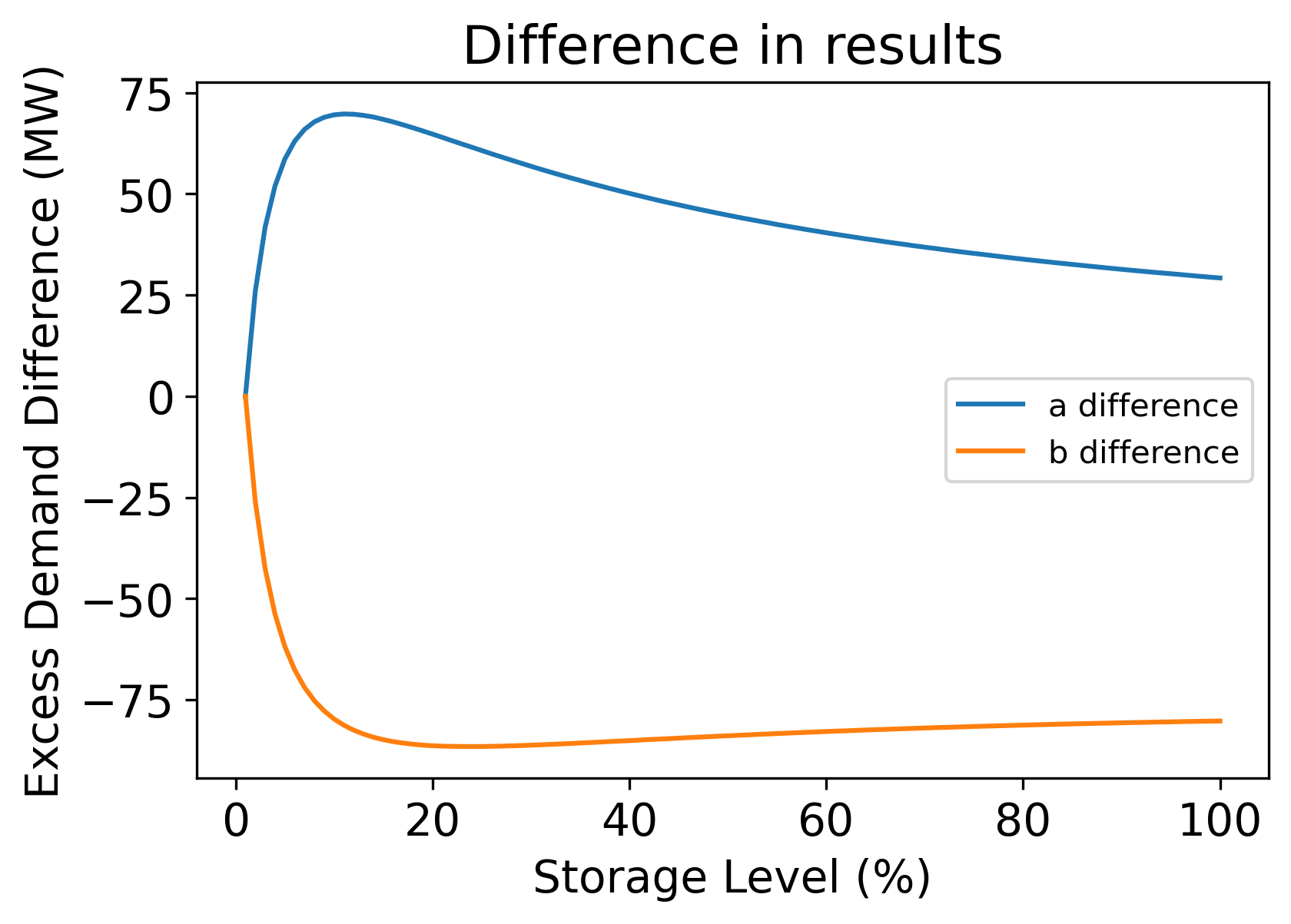} }
	\caption{Comparison of results between explicit calculation and numerical method}
\end{figure}

It can be seen that the results of the numerical method closely follow the true values of the optimal control points, with $a$ being at most $69.76$MW, away from the true solution, and $b$ at most $-86.56$MW. These errors can be reduced by using a finer grid, as in this case $Z = \left[ 1\%, 2\%, ... 100\% \right]$ with $\Delta z = 1\%$.

\newpage
\subsection{Variation with Temperature}
Of particular interest for power storage facility operators is how to operate the facility throughout the year. In particular, distinctions can be made between Summer, Winter, and the shoulder months comprising Autumn and Spring. The shoulder months are of particular interest, as temperature can fluctuate significantly over the course of a week leading to increased uncertainty in the future price of electricity. This seasonal grouping can be considered a factor of temperature, with more electricity required to deal with the cold during the winter months. The effects of temperature could be considered using a 2 dimensional stochastic process underlying the data, however the analysis required to obtain applicable results for such a process is beyond the scope of this paper. 

Instead, the effect of temperature, $T$, can be considered through a translation on the payoff functions $E, F$ to account for the increase in demand. As the temperature decreases, excess demand increases and leads to more expensive forms of electricity generation being required, increasing the price of electricity. In particular, there is a larger growth in demand for higher values than lower, meaning some scaling is necessary to account for this. Furthermore the increase in demand is more exponential than linear, and therefore we consider payoff functions of the following form:
\begin{align*}
    F(x,T) &= 0.001\left(x+(20-T)^2\left(55 + \frac{x-min(x)}{1500}\right)\right)+20 \\
    E(x,T) &= 0.9 \times \left( 0.001\left(x+(20-T)^2\left(55 + \frac{x-min(x)}{1500}\right)\right)+20 \right)
\end{align*}

where $T \in [5, 20]$ is the average temperature over a month, and $min(x) = -40,000$ is the minimum excess demand value in the range. Whilst daily temperature may be a better predictor for excess demand for any particular day, average monthly temperature provides a simpler relationship between demand and temperature, and is a suitable indicator for the relationship between temperature and excess demand.

\begin{figure}[h!]
    \centering
    \includegraphics[width=0.75\textwidth]{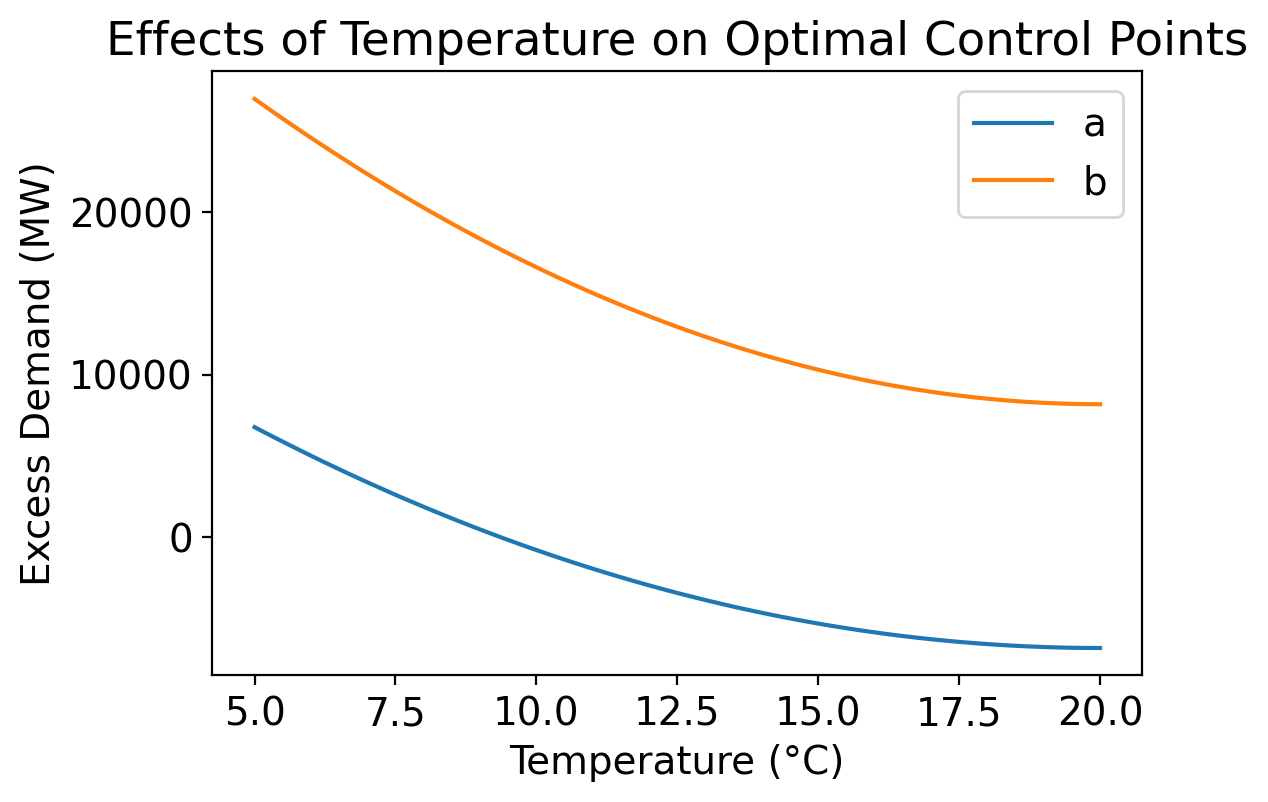}
    \caption{The effect of temperature on the optimal control points.}
\end{figure}

From Figure 4, it can be seen both optimal control points have an inverse relationship with temperature, increasing exponentially as the total level of excess demand increases. This is expected as the optimal control points move to remain in appropriate levels for the increasing mean level of excess demand. Furthermore it can be noted that the optimal control points increase at different rates with the gap between them increasing with the average level of excess demand, also expected as the increase in demand is larger at higher levels of excess demand at lower levels.

\subsection{Composition of Different Variables}

Assume that the variable payoff functions are functions of excess demand $x$, factors which are separable from excess demand, $\mathbf{y} = (y_1, y_2, ....)$, and factors which cannot, $\mathbf{u} = (u_1, u_2, ....)$, with all variables being independent of each other. The variation with the factor $y_j$ on the payoff functions $F, E$ can be represented by the functions $Y_{j,f}(y_j), Y_{j,e}(y_j)$ respectively, resulting in:
\begin{align*}
    F(x, \mathbf{y}, \mathbf{u}) &= \prod_j{Y_{j,f}(y_j)}F_x(x, \mathbf{u}) \\
    E(x, \mathbf{y}, \mathbf{u}) &= \prod_j{Y_{j,e}(y_j)}E_x(x, \mathbf{u}).
\end{align*}

\newpage
To simplify the problem, whilst considering the variable $y_k$ the rest are held constant. We are interested in the variation of the optimal control points with regards to the variable $y_k$, meaning that the optimal control strategy can be represented by $(a(y_k(i_k),\mathbf{y}^*, \mathbf{u})), b(y_k(i_k + 1),\mathbf{y}^*, \mathbf{u}))$ where $\mathbf{y}^*$ denotes the vector $\mathbf{y}$ without the variable $y_k$. For ease of notation, we now define:
\begin{align*}
    &a(y_k(i_k),\mathbf{y}^*, \mathbf{u}))  = a_i \\
    &b(y_k(i_k + 1),\mathbf{y}^*, \mathbf{u}) = b_{i+1}
\end{align*}

Thus the system of equations (5-6) becomes:
\begin{align}
    \prod_j{Y_{j,f}(y_j(i_j))}\frac{d}{dx}\left( \frac{F_x(x, \mathbf{u})}{\psi(x)} \right)_{x=a_i} \frac{\psi^2(a_i)}{\mathcal{W}(a_i)} &= \prod_j{Y_{j,e}(y_j(i_j+1))}\frac{d}{dx}\left( \frac{E_{x}(x, \mathbf{u})}{\psi(x)} \right)_{x=b_{i+1}} \frac{\psi^2(b_{i+1})}{\mathcal{W}(b_{i+1})} \\
    \prod_j{Y_{j,f}(y_j(i_j))}\frac{d}{dx}\left( \frac{F_x(x, \mathbf{u})}{\phi(x)} \right)_{x=a_i} \frac{\phi^2(a_i)}{\mathcal{W}(a_i)} &= \prod_j{Y_{j,e}(y_j(i_j+1))}\frac{d}{dx}\left( \frac{E_{x}(x, \mathbf{u})}{\phi(x)} \right)_{x=b_{i+1}} \frac{\phi^2(b_{i+1})}{\mathcal{W}(b_{i+1})}
\end{align}

For ease of notation we now make the following definitions:
\begin{align*}
    q_{\psi, f}(x, \mathbf{y}^*, \mathbf{u}) &= \prod_{j \ne k}{Y_{j,f}(y_j(i_j))}\frac{d}{dx}\left( \frac{F_x(x, \mathbf{u})}{\psi(x)} \right) \frac{\psi^2(x)}{\mathcal{W}(x)} = \frac{d}{dx}\left( \frac{F_x(x, \mathbf{y}^*, \mathbf{u})}{\psi(x)} \right) \frac{\psi^2(x)}{\mathcal{W}(x)}\\
    q_{\psi, e}(x, \mathbf{y}^*, \mathbf{u}) &= \prod_{j \ne k}{Y_{j,e}(y_j(i_j))}\frac{d}{dx}\left( \frac{E_x(x, \mathbf{u})}{\psi(x)} \right) \frac{\psi^2(x)}{\mathcal{W}(x)} = \frac{d}{dx}\left( \frac{F_x(x, \mathbf{y}^*, \mathbf{u})}{\psi(x)} \right) \frac{\psi^2(x)}{\mathcal{W}(x)}\\
    q_{\phi, f}(x, \mathbf{y}^*, \mathbf{u}) &= \prod_{j \ne k}{Y_{j,f}(y_j(i_j))}\frac{d}{dx}\left( \frac{F_x(x, \mathbf{u})}{\phi(x)} \right) \frac{\phi^2(x)}{\mathcal{W}(x)} = \frac{d}{dx}\left( \frac{F_x(x, \mathbf{y}^*, \mathbf{u})}{\psi(x)} \right) \frac{\psi^2(x)}{\mathcal{W}(x)}\\
    q_{\phi, e}(x, \mathbf{y}^*, \mathbf{u}) &= \prod_{j \ne k}{Y_{j,e}(y_j(i_j))}\frac{d}{dx}\left( \frac{E_x(x, \mathbf{u})}{\phi(x)} \right) \frac{\phi^2(x)}{\mathcal{W}(x)} = \frac{d}{dx}\left( \frac{F_x(x, \mathbf{y}^*, \mathbf{u})}{\psi(x)} \right) \frac{\psi^2(x)}{\mathcal{W}(x)}
\end{align*}

such that we can re-write the system of equations as:
\begin{align*}
    Y_{k,f}(y_k(i_k))q_{\psi, f}(a_i, \mathbf{y}^*, \mathbf{u}) &= Y_{k,e}(y_k(i_k+1))q_{\psi, e}(b_{i+1}, \mathbf{y}^*, \mathbf{u}) \\
    Y_{k,f}(y_k(i_k))q_{\phi, f}(a_i, \mathbf{y}^*, \mathbf{u}) &= Y_{j,e}(y_k(i_k + 1))q_{\phi, e}(b_{i+1}, \mathbf{y}^*, \mathbf{u})
\end{align*}

Due to the assumed independence of $y_k$ with $x, \mathbf{y}^*, \mathbf{u}$, we can proceed as previously to find the following:
\begin{align*}
    \frac{da}{dy_k} &= \left[ \frac{Y'_{k,e}(y_k)}{Y_{k,e}(y_k)} - \frac{Y'_{k,f}(y_k)}{Y_{k,f}(y_k)} \right] \left( \frac{q_{\psi,f}(a(x, \mathbf{y}^*, \mathbf{u})q'_{\phi,e}(b(x, \mathbf{y}^*, \mathbf{u}) - q'_{\psi,e}(b(x, \mathbf{y}^*, \mathbf{u})q_{\phi,f}(a(x, \mathbf{y}^*, \mathbf{u})}{q'_{\psi,f}(a(x, \mathbf{y}^*, \mathbf{u})q'_{\phi,e}(b(x, \mathbf{y}^*, \mathbf{u}) - q'_{\psi,e}(b(x, \mathbf{y}^*, \mathbf{u})q'_{\phi,f}(a(x, \mathbf{y}^*, \mathbf{u})} \right) \\
    &= \frac{d}{dy_k}\left( \frac{Y_{k,e}(y_k)}{Y_{k,f}(y_k)} \right)\left( \frac{q_{\psi,e}(b(x, \mathbf{y}^*, \mathbf{u})q'_{\phi,e}(b(x, \mathbf{y}^*, \mathbf{u}) - q'_{\psi,e}(b(x, \mathbf{y}^*, \mathbf{u})q_{\phi,e}(b(x, \mathbf{y}^*, \mathbf{u})}{q'_{\psi,f}(a(x, \mathbf{y}^*, \mathbf{u})q'_{\phi,e}(b(x, \mathbf{y}^*, \mathbf{u}) - q'_{\psi,e}(b(x, \mathbf{y}^*, \mathbf{u})q'_{\phi,f}(a(x, \mathbf{y}^*, \mathbf{u})} \right) \numberthis
\end{align*}
\begin{align*}
    \frac{db}{dy_k} &= \left[ \frac{Y'_{k,e}(y_k)}{Y_{k,e}(y_k)} - \frac{Y'_{k,f}(y_k)}{Y_{k,f}(y_k)} \right] \left( \frac{q_{\psi,e}(b(x, \mathbf{y}^*, \mathbf{u})q'_{\phi,f}(a(x, \mathbf{y}^*, \mathbf{u}) - q'_{\psi,f}(a(x, \mathbf{y}^*, \mathbf{u})q_{\phi,e}(b(x, \mathbf{y}^*, \mathbf{u})}{q'_{\psi,f}(a(x, \mathbf{y}^*, \mathbf{u})q'_{\phi,e}(b(x, \mathbf{y}^*, \mathbf{u}) - q'_{\psi,e}(b(x, \mathbf{y}^*, \mathbf{u})q'_{\phi,f}(a(x, \mathbf{y}^*, \mathbf{u})} \right) \\
    &= \frac{d}{dy_k}\left( \frac{Y_{k,f}(y_k)}{Y_{k,e}(y_k)} \right)\left( \frac{q_{\psi,f}(a(x, \mathbf{y}^*, \mathbf{u})q'_{\phi,f}(a(x, \mathbf{y}^*, \mathbf{u}) - q'_{\psi,f}(a(x, \mathbf{y}^*, \mathbf{u})q_{\phi,f}(a(x, \mathbf{y}^*, \mathbf{u})}{q'_{\psi,f}(a(x, \mathbf{y}^*, \mathbf{u})q'_{\phi,e}(b(x, \mathbf{y}^*, \mathbf{u}) - q'_{\psi,e}(b(x, \mathbf{y}^*, \mathbf{u})q'_{\phi,f}(a(x, \mathbf{y}^*, \mathbf{u})} \right). \numberthis
\end{align*}
By varying one factor at a time, it is possible to identify an n-dimensional manifold which defines the optimal control strategy for all variations in the factors $\mathbf{y}, \mathbf{u}$. As an example, consider the effects of storage level and temperature considered previously. Varying storage level $z$ for all temperatures $T$ allows for optimal control of the power storage facility for all storage levels throughout the year. If we consider the functions seen previously for both temperature and storage levels, then the payoff functions can be defined as below:
\begin{align*}
    F(x, z, T) &= \left( 0.001\left(x+(20-T)^2\left(55 + \frac{x-min(x)}{1500}\right)\right)+20 \right) \left(1 + 2\frac{z}{z_n} \right) \\
    E(x, z, T) &= 0.9\left( 0.001\left(x+(20-T)^2\left(55 + \frac{x-min(x)}{1500}\right)\right)+20 \right)
\end{align*}

Surfaces can be found for the optimal control points $a, b$ as seen in Figure 5 for these payoff functions.

\begin{figure}[h!]
    \centering
    \includegraphics[width=0.75\textwidth]{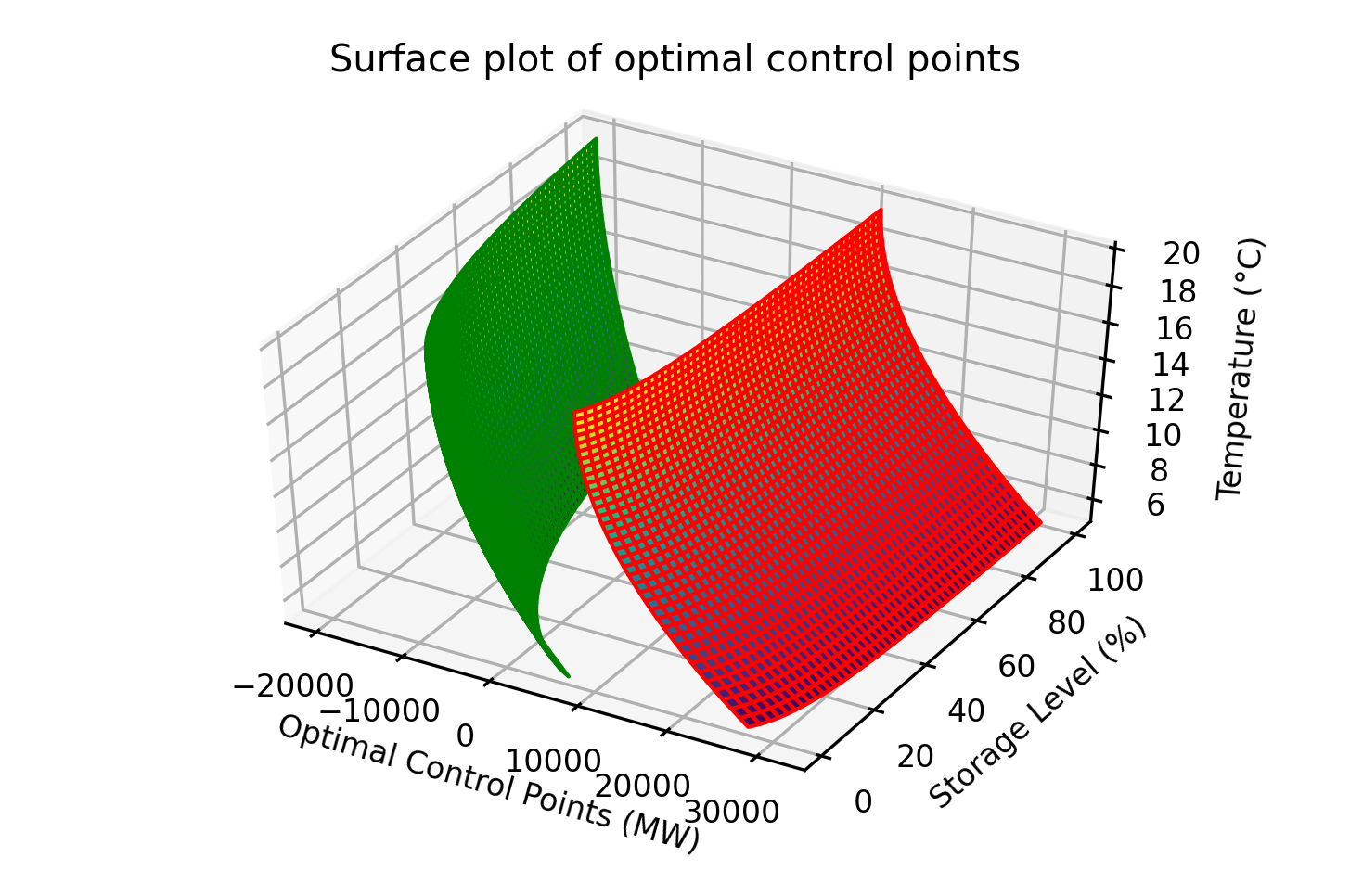}
    \caption{Surface plot showing the interaction between storage level and Temperature.}
\end{figure}

The trends seen previously for both storage level and temperature can be observed in these surface plots, though it can be seen that there is interaction between the variables. In particular each optimal control point is affected more by one variable than the other, as can be seen by comparing the end cases for the facility as in Figures 6 \& 7.

\newpage
\begin{figure}[h!]
    \centering
    \includegraphics[width=0.9\textwidth]{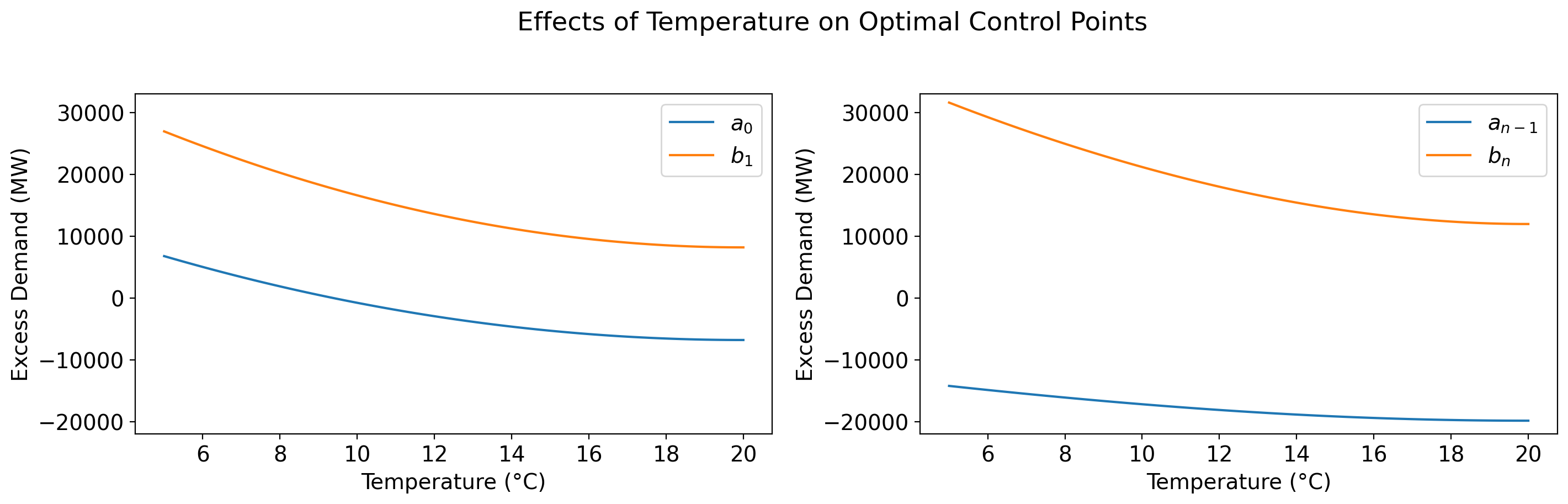}
    \caption{Comparison of the effect of temperature on $(a_0, b_1)$ and $(a_{n-1},b_n)$.}
\end{figure}

Whilst the exponential trends can still be seen on $(a_{n-1}, b_n)$, the variation in $a_{n-1}$ is noticeably less than for $a_0$. This is likely due to the large effect that storage level has on $a$, whilst $b$ varies with storage level up to approximately $20\%$ level before plateauing.

\begin{figure}[h!]
    \centering
    \includegraphics[width=0.9\textwidth]{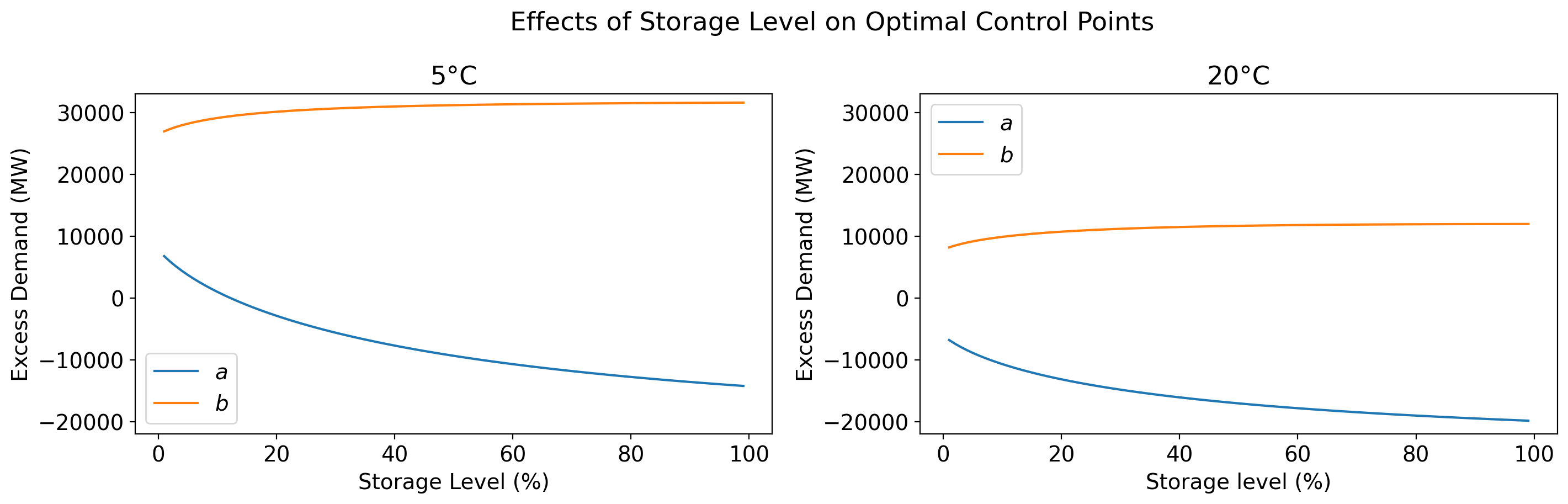}
    \caption{Comparison of the effect of Storage Level at different Temperatures.}
\end{figure}

By comparing the effects of storage level at $5$°C and $20$°C, it can be seen that temperature has a larger effect on the upper control point $b$ and whilst $a$ is initially effected at lower levels of storage it is forced to more extreme levels of excess demand by the increasing cost of filling the facility.

\section{Discussion}

The work of \cite{EG} and \cite{EG-J} developed an approach for identifying the optimal control points for a facility directly from excess demand data, but was naive in identifying only one pair of optimal control points for a facility without consideration for the variation of a range of important factors. The developments within this paper allow for an optimal control strategy to be identified that is tailored to a specific facility, with a consideration of potential other factors. Such developments could also be applied to other scenarios outside of power storage facilities, such as those in discussed in the context of ecological and environmental management in \cite{Y-Y}

\newpage
An important simplification made in this paper was the assumption of a linear payoff function, which can be complicated by the inclusion of price jumps seen in the MOC. Future work could consider how best to handle the complications introduced by working with a more complicated payoff function, and in particular how to handle scenarios in which there are multiple identified control strategies which satisfy equations (5-6), (13-14) or (17-18).

\textbf{Acknowledgement of Funding:}

FO'B was supported by the EPSRC Centre for Doctoral Training in Mathematical Modelling, Analysis and Computation (MAC-MIGS) funded by the UK Engineering and Physical Sciences Research Council (grant EP/S023291/1), Heriot-Watt University and the University of Edinburgh.

\newpage

\section{Appendix A}
We begin by applying the principal of smooth fit to our end cases to determine the values of the optimal control points, starting with $z_0$:
\begin{align*}
A_1\phi(a_0) + B_1\psi(a_0) - F_0(a_0) &= A_0\phi(a_0) \\
A_1\phi'(a_0) + B_1\psi'(a_0) - F_0'(a_0) &= A_0\phi'(a_0)
\end{align*}
gathering similar terms gives:
\begin{align}
&\left[ A_1 - A_0 \right] \phi(a_0) + B_1\psi(a_0)  = F_0(a_0) \\
&\left[ A_1 - A_0 \right] \phi'(a_0) + B_1\psi'(a_0) = F_0'(a_0) 
\end{align}
rearranging (21) for $\left[ A_1 - A_0 \right]$:
$$\left[ A_1 - A_0 \right] = \frac{F_0(a_0) - B_1\psi(a_0)}{\phi(a_0)} $$
substituting into (22) and rearranging:
\begin{align*}
    &\frac{F_0(a_0) - B_1\psi(a_0)}{\phi(a_0)} \phi'(a_0) + B_1\psi'(a_0) = F_0'(a_0) \\
    \implies &F_0'(a_0)\phi(a_0) - B_1\psi'(a_0)\phi(a_0)  =  F_0(a_0)\phi'(a_0) - B_1\psi(a_0)\phi'(a_0) \\
    \implies &B_1 \left( \psi'(a_0)\phi(a_0) - \psi(a_0)\phi'(a_0) \right) = F'_0(a_0) \phi(a_0) - F_0(a_0)\phi'(a_0)  \\
    \implies &B_1 \mathcal{W}(a_0) = \frac{d}{dx}\left( \frac{F_0(x)}{\phi(x)} \right)_{x=a_0} \phi^2(a_0) \\
    \implies &B_1 = \frac{d}{dx}\left( \frac{F_0(x)}{\phi(x)} \right)_{x=a_0} \frac{\phi^2(a_0)}{\mathcal{W}(a_0)}.
\end{align*} 

Alternatively (21) can be rearranged for $B_1$:
$$ B_1  = \frac{F_0(a_0) - \left[ A_1 - A_0 \right] \phi(a_0)}{\psi(a_0)}$$
substituting into (22) and rearranging:
\begin{align*}
    &\left[ A_1 - A_0 \right] \phi'(a_0) + \frac{F_0(a_0) - \left[ A_1 - A_0 \right] \phi(a_0)}{\psi(a_0)}\psi'(a_0) = F_0'(a_0) \\
    \implies &\left[ A_1 - A_0 \right] \psi(a_0)\phi'(a_0) + F_0(a_0)\psi'(a_0) - \left[ A_1 - A_0 \right] \psi'(a_0)\phi(a_0) = F_0'(a_0)\psi(a_0) \\
    \implies &\left[ A_1 - A_0 \right] \left( \psi(a_0)\phi'(a_0) - \psi'(a_0)\phi(a_0) \right) = F_0'(a_0)\psi(a_0) - F_0(a_0)\psi'(a_0)  \\
    \implies &-\left[ A_1 - A_0 \right] \mathcal{W}(a_0) = \frac{d}{dx} \left( \frac{F_0(x)}{\psi(x)} \right)_{x = a_0} \psi^2(a_0) \\
    \implies &-\left[ A_1 - A_0 \right] = \frac{d}{dx} \left( \frac{F_0(x)}{\psi(x)} \right)_{x = a_0} \frac{\psi^2(a_0)}{\mathcal{W}(a_0)}
\end{align*}

\newpage

We can similarly consider the end case at $z_n$:
\begin{align*}
B_n\psi(b_n)  &= A_{n-1}\phi(b_n) + B_{n-1}\psi(b_n) + E_{n}(b_n) \\
B_n\psi'(b_n)  &= A_{n-1}\phi'(b_n) + B_{n-1}\psi'(b_n) + E_{n}'(b_n)
\end{align*}
gathering similar terms:
\begin{align}
&\left[ B_n - B_{n-1} \right] \psi(b_n) - A_{n-1}\phi(b_n) =  E_{n}(b_n) \\
&\left[ B_n - B_{n-1} \right] \psi'(b_n) - A_{n-1}\phi'(b_n) =  E_{n}'(b_n)
\end{align}
rearranging (23) for $\left[ B_n - B_{n-1} \right]$:
$$\left[ B_n - B_{n-1} \right]   =  \frac{E_{n}(b_n) + A_{n-1}\phi(b_n)}{\psi(b_n)} $$
substituting into (24) and rearranging:
\begin{align*}
    &\frac{E_{n}(b_n) + A_{n-1}\phi(b_n)}{\psi(b_n)}\psi'(b_n) - A_{n-1}\phi'(b_n) =  E_{n}'(b_n) \\
    \implies &E_n(b_n)\psi'(b_n) + A_{n-1}\psi'(b_n)\phi(b_n) = E_n'(b_n)\psi(b_n) + A_{n-1}\psi(b_n)\phi'(b_n) \\
    \implies & A_{n-1}\left( \psi'(b_n)\phi(b_n) - \psi(b_n)\phi'(b_n) \right) = E_n'(b_n)\psi(b_n) - E_n(b_n)\psi'(b_n) \\
    \implies &A_{n-1}\mathcal{W}(b_n) = \frac{d}{dx}\left(\frac{E_n(x)}{\psi(x)}\right)_{x=b_n}\psi^2(b_n)\\
    \implies &A_{n-1} = \frac{d}{dx}\left(\frac{E_n(x)}{\psi(x)}\right)_{x=b_n}\frac{\psi^2(b_n)}{\mathcal{W}(b_n)}
\end{align*}
Alternatively (23) can be rearranged for $A_{n-1}$:
$$ A_{n-1} = \frac{\left[ B_n - B_{n-1} \right] \psi(b_n) -  E_{n}(b_n)}{\phi(b_n)}$$
substituting into (24) and rearranging:
\begin{align*}
    &\left[ B_n - B_{n-1} \right] \psi'(b_n) - \frac{\left[ B_n - B_{n-1} \right] \psi(b_n) -  E_{n}(b_n)}{\phi(b_n)}\phi'(b_n) =  E_{n}'(b_n) \\
    \implies & \left[ B_n - B_{n-1} \right] \psi'(b_n)\phi(b_n) - \left[ B_n - B_{n-1} \right] \psi(b_n)\phi'(b_n) -  E_{n}(b_n)\phi'(b_n) =  E_{n}'(b_n)\phi(b_n) \\
    \implies & \left[ B_n - B_{n-1} \right] \left( \psi'(b_n)\phi(b_n) - \psi(b_n)\phi'(b_n) \right)  =  E_{n}'(b_n)\phi(b_n) + E_{n}(b_n)\phi'(b_n) \\
    \implies & \left[ B_n - B_{n-1} \right] \mathcal{W}(b_n)  =  \frac{d}{dx} \left( \frac{E_n(x)}{\phi(x)} \right)_{x=b_n} \phi^2(b_n) \\
    \implies & \left[ B_n - B_{n-1} \right] =  \frac{d}{dx} \left( \frac{E_n(x)}{\phi(x)} \right)_{x=b_n} \frac{\phi^2(b_n)}{\mathcal{W}(b_n)}
\end{align*}

\newpage

We next consider the boundary values for the value function for some storage value $z_i, 0 < i < n$:
\begin{align*}
    w(x,z_i) = 
          \begin{dcases}
            A_{i+1}\phi(x) + B_{i+1}\psi(x) - F(x,z_i), &\text{ if } x \in (\alpha, a_i] \\
            A_i\phi(x) + B_i\psi(x), &\text{ if } x \in (a_i, b_i) \\
            A_{i-1}\phi(x) + B_{i-1}\psi(x) + E(x,z_i), &\text{ if } x \in [b_i, \beta) 
            \end{dcases}
\end{align*}

Consider the lower boundary $a_i$, by applying the principle of smooth fit we have:
\begin{align*}
A_{i+1}\phi(a_i) + B_{i+1}\psi(a_i) - F_i(a_i) &= A_i\phi(a_i) + B_i\psi(a_i) \\
A_{i+1}\phi'(a_i) + B_{i+1}\psi'(a_i) - F'_i(a_i) &= A_i\phi'(a_i) + B_i\psi'(a_i)
\end{align*}
grouping similar terms:
\begin{align}
&\left[ A_{i+1} -  A_i \right] \phi(a_i) + \left[ B_{i+1} - B_i \right] \psi(a_i)  =  F_i(a_i) \\
&\left[ A_{i+1} -  A_i \right] \phi'(a_i) + \left[ B_{i+1} - B_i \right] \psi'(a_i)  =  F_i'(a_i)
\end{align}
rearranging (25) for $\left[ A_{i+1} -  A_i \right]$:
$$\left[ A_{i+1} -  A_i \right]  =  \frac{F_i(a_i) - \left[ B_{i+1} - B_i \right] \psi(a_i)}{\phi(a_i)}$$
substituting into (26) and rearranging
\begin{align*}
    &\frac{F_i(a_i) - \left[ B_{i+1} - B_i \right] \psi(a_i)}{\phi(a_i)} \phi'(a_i) + \left[ B_{i+1} - B_i \right] \psi'(a_i)  =  F_i'(a_i) \\
    \implies &F_i(a_i)\phi'(a_i) - \left[ B_{i+1} - B_i \right] \psi(a_i) \phi'(a_i) + \left[ B_{i+1} - B_i \right] \psi'(a_i)\phi(a_i)  =  F_i'(a_i)\phi(a_i) \\
    \implies &\left[ B_{i+1} - B_i \right] \left(\psi'(a_i)\phi(a_i) - \psi(a_i) \phi'(a_i) \right) =  F_i'(a_i)\phi(a_i) - F_i(a_i)\phi'(a_i) \\
    \implies & \left[ B_{i+1} - B_i \right]\mathcal{W}(a_i) = \frac{d}{dx}\left( \frac{F_i(x)}{\phi(x)} \right)_{x=a_i} \phi^2(a_i) \\
    \implies & \left[ B_{i+1} - B_i \right] = \frac{d}{dx}\left( \frac{F_i(x)}{\phi(x)} \right)_{x=a_i} \frac{\phi^2(a_i)}{\mathcal{W}(a_i)}.
\end{align*}

Alternatively we can rearrange (25) for $\left[ B_{i+1} - B_i \right]$:
$$ \left[ B_{i+1} - B_i \right]   =  \frac{F_i(a_i) - \left[ A_{i+1} -  A_i \right] \phi(a_i)}{\psi(a_i)}$$
and in a similar manner to before, substitute into (26) to find:
\begin{align*}
    &\left[ A_{i+1} -  A_i \right] \phi'(a_i) +  \frac{F_i(a_i) - \left[ A_{i+1} -  A_i \right] \phi(a_i)}{\psi(a_i)}\psi'(a_i)  =  F_i'(a_i)\\
    \implies &\left[ A_{i+1} -  A_i \right] \psi(a_i)\phi'(a_i) +  F_i(a_i)\psi'(a_i) - \left[ A_{i+1} -  A_i \right] \phi(a_i)\psi'(a_i)  =  F_i'(a_i)\psi(a_i) \\
    \implies &\left[ A_{i+1} -  A_i \right] \left( \psi(a_i)\phi'(a_i) - \phi(a_i)\psi'(a_i)\right)  =  F_i'(a_i)\psi(a_i) - F_i(a_i)\psi'(a_i) \\
    \implies &-\left[ A_{i+1} -  A_i \right] \mathcal{W}(a_i)  =  \frac{d}{dx}\left( \frac{F_i(x)}{\psi(x)} \right)_{x=a_i}\psi^2(a_i) \\
    \implies &-\left[ A_{i+1} -  A_i \right] = \frac{d}{dx}\left( \frac{F_i(x)}{\psi(x)} \right)_{x=a_i}\frac{\psi^2(a_i)}{\mathcal{W}(a_i)}
\end{align*}

\newpage
Similarly consider the upper boundary $b_{i}$, where we have:
\begin{align*}
    A_i\phi(b_i) + B_i\psi(b_i) &= A_{i-1}\phi(b_i) + B_{i-1}\psi(b_i) + E_i(b_i) \\
    A_i\phi'(b_i) + B_i\psi'(b_i) &= A_{i-1}\phi'(b_i) + B_{i-1}\psi'(b_i) + E'_i(b_i)
\end{align*}
grouping similar terms:
\begin{align}
    &\left[ A_i - A_{i-1} \right] \phi(b_i) + \left[ B_i - B_{i-1} \right] \psi(b_i) =  E_i(b_i) \\
    &\left[ A_i - A_{i-1} \right] \phi'(b_i) + \left[ B_i - B_{i-1} \right] \psi'(b_i) =  E_i'(b_i)
\end{align}
rearranging (27) for $\left[ A_i - A_{i-1} \right]$:
$$\left[ A_i - A_{i-1} \right]   =  \frac{E_i(b_i) - \left[ B_i - B_{i-1} \right] \psi(b_i)}{\phi(b_i)}$$
substituting into (28) and rearranging:
\begin{align*}
    &\frac{E_i(b_i) - \left[ B_i - B_{i-1} \right] \psi(b_i)}{\phi(b_i)} \phi'(b_i) + \left[ B_i - B_{i-1} \right] \psi'(b_i) =  E_i'(b_i) \\
    \implies &E_i(b_i)\phi'(b_i) - \left[ B_i - B_{i-1} \right] \psi(b_i) \phi'(b_i) + \left[ B_i - B_{i-1} \right] \psi'(b_i)\phi(b_i) =  E_i'(b_i)\phi(b_i)\\
    \implies & \left[ B_i - B_{i-1} \right] \left( \psi'(b_i)\phi(b_i) -  \psi(b_i) \phi'(b_i) \right) =  E_i'(b_i)\phi(b_i) - E_i(b_i)\phi'(b_i) \\
    \implies & \left[ B_i - B_{i-1} \right] \mathcal{W}(b_i) =  \frac{d}{dx} \left( \frac{E_i(x)}{\phi(x)} \right)_{x=b_i} \phi^2(b_i) \\
    \implies & \left[ B_i - B_{i-1} \right]  =  \frac{d}{dx} \left( \frac{E_i(x)}{\phi(x)} \right)_{x=b_i} \frac{\phi^2(b_i)}{\mathcal{W}(b_i)}.
\end{align*}

Alternatively we can rearrange (27) for $\left[ B_i - B_{i-1} \right]$:
$$\left[ B_i - B_{i-1} \right]  = \frac{E_i(b_i) - \left[ A_i - A_{i-1} \right] \phi(b_i)}{\psi(b_i)}$$
and in a similar manner to before, substitute into (28) to find:
\begin{align*}
    &\left[ A_i - A_{i-1} \right] \phi'(b_i) + \frac{E_i(b_i) - \left[ A_i - A_{i-1} \right] \phi(b_i)}{\psi(b_i)} \psi'(b_i) =  E_i'(b_i) \\
    \implies &\left[ A_i - A_{i-1} \right] \phi'(b_i)\psi(b_i) + E_i(b_i)\psi'(b_i) - \left[ A_i - A_{i-1} \right] \phi(b_i) \psi'(b_i) =  E_i'(b_i)\psi(b_i) \\
    \implies &\left[ A_i - A_{i-1} \right] \left( \phi'(b_i)\psi(b_i) - \phi(b_i) \psi'(b_i) \right) =  E_i'(b_i)\psi(b_i) - E_i(b_i)\psi'(b_i) \\
    \implies &-\left[ A_i - A_{i-1} \right] \mathcal{W}(b_i) = \frac{d}{dx} \left( \frac{E_i(x)}{\psi(x))} \right)_{x=b_i} \psi^2(b_i) \\
    \implies &-\left[ A_i - A_{i-1} \right]  = \frac{d}{dx} \left( \frac{E_i(x)}{\psi(x))} \right)_{x=b_i} \frac{\psi^2(b_i)}{\mathcal{W}(b_i)}
\end{align*}

\newpage

Thus we have the following equations:
\begin{align}
    -&\left[ A_{i+1} -  A_i \right] = \frac{d}{dx}\left( \frac{F_i(x)}{\psi(x)} \right)_{x=a_i}\frac{\psi^2(a_i)}{\mathcal{W}(a_i)}\\
    -&\left[ A_i - A_{i-1} \right]  = \frac{d}{dx} \left( \frac{E_i(x)}{\psi(x))} \right)_{x=b_i} \frac{\psi^2(b_i)}{\mathcal{W}(b_i)}\\
    &\left[ B_{i+1} - B_i \right] = \frac{d}{dx}\left( \frac{F_i(x)}{\phi(x)} \right)_{x=a_i} \frac{\phi^2(a_i)}{\mathcal{W}(a_i)}\\
    &\left[ B_i - B_{i-1} \right]  =  \frac{d}{dx} \left( \frac{E_i(x)}{\phi(x)} \right)_{x=b_i} \frac{\phi^2(b_i)}{\mathcal{W}(b_i)}
\end{align}
for $1 \le i \le {n-1}$. setting $i = 1$ combined with our knowledge at the end case $i = 0$ we can see that:
\begin{align*}
- \left[ A_1 - A_0 \right] &= \frac{d}{dx} \left( \frac{F_0(x)}{\psi(x)} \right)_{x = a_0} \frac{\psi^2(a_0)}{\mathcal{W}(a_0)} = \frac{d}{dx} \left( \frac{E_1(x)}{\psi(x))} \right)_{x=b_1} \frac{\psi^2(b_1)}{\mathcal{W}(b_1)} \\
B_1 &= \frac{d}{dx}\left( \frac{F_0(x)}{\phi(x)} \right)_{x=a_0} \frac{\phi^2(a_0)}{\mathcal{W}(a_0)} = \frac{d}{dx} \left( \frac{E_1(x)}{\phi(x))} \right)_{x=b_1} \frac{\phi^2(b_1)}{\mathcal{W}(b_1)}
\end{align*}

similarly setting $i = n-1$ combined with our knowledge at the end case $i = n$:
\begin{align*}
    A_{n-1} &= \frac{d}{dx} \left( \frac{F_{n-1}(x)}{\psi(x)} \right)_{x=a_{n-1}} \frac{\psi^2(a_{n-1})}{\mathcal{W}(a_{n-1})}  = \frac{d}{dx}\left(\frac{E_n(x)}{\psi(x)}\right)_{x=b_n}\frac{\psi^2(b_n)}{\mathcal{W}(b_n)}\\
    \left[ B_n - B_{n-1} \right] &= \frac{d}{dx} \left( \frac{F_{n-1}(x)}{\phi(x)} \right)_{x=a_{n-1}} \frac{\phi^2(a_{n-1})}{\mathcal{W}(a_{n-1})} = \frac{d}{dx} \left( \frac{E_n(x)}{\phi(x)} \right)_{x=b_n} \frac{\phi^2(b_n)}{\mathcal{W}(b_n)}.
\end{align*}

If we were to consider the control point $b_{i+1}$ instead for equations (29-32), then for $1 \le i \le n-2$ we would have the system of equations given by:
\begin{align}
    -\left[ A_{i+1}-A_{i} \right] &= \frac{d}{dx}\left( \frac{F_i(x)}{\psi(x)} \right)_{x = a_i} \frac{\psi^2(a_i)}{\mathcal{W}(a_i)} = \frac{d}{dx}\left( \frac{E_{i+1}(x)}{\psi(x)} \right)_{x = b_{i+1}} \frac{\psi^2(b_{i+1})}{\mathcal{W}(b_{i+1})} \\
    \left[ B_{i+1}-B_{i} \right] &= \frac{d}{dx}\left( \frac{F_i(x)}{\phi(x)} \right)_{x = a_i} \frac{\phi^2(a_i)}{\mathcal{W}(a_i)} = \frac{d}{dx}\left( \frac{E_{i+1}(x)}{\phi(x)} \right)_{x=b_{i+1}} \frac{\phi^2(b_{i+1})}{\mathcal{W}(b_{i+1})}
\end{align}

which combined with our equations at the end cases $z_0, z_n$ and the knowledge that $A_n = B_0 = 0$ can be expanded to $0 \le i \le i-1$. Thus it is possible to identify the optimal control points for the full operating range of the power storage facility.

\newpage

\section{Appendix B}

Assume that for some value $1 \le i \le n-2$ the optimal control points $(a(z_i), b(z_{i+1}))$ solving 
\begin{align}
    & Z_f(z_{i})q_{\psi, f}(a_{i}) = Z_e(z_{i+1})q_{\psi, e}(b_{i+1}) \\
    & Z_f(z_{i})q_{\phi, f}(a_{i}) = Z_e(z_{i+1})q_{\phi, e}(b_{i+1})
\end{align}
are known, and we aim to identify the optimal control points $(a(z_{i+1}), b(z_{i+2}))$ which satisfy:
\begin{align}
    & Z_f(z_{i+1})q_{\psi, f}(a_{i+1}) = Z_e(z_{i+2})q_{\psi, e}(b_{i+2}) \\
    & Z_f(z_{i+1})q_{\phi, f}(a_{i+1}) = Z_e(z_{i+2})q_{\phi, e}(b_{i+2})
\end{align}

Assume that $z_{i+1}-z_i = z_{i+2}-z_{i+1} = \Delta z$ corresponds to $a(z_{i+1})-a(z_{i}) = \Delta a, b(z_{i+2}) - b(z_{i+1}) = \Delta b$, then we can rewrite equations (37-38) as:
\begin{align}
    & Z_f(z_i + \Delta z)q_{\psi, f}(a(z_{i}) + \Delta a) = Z_e(z_{i+1} + \Delta z)q_{\psi, e}(b(z_{i+1}) + \Delta b) \\
    & Z_f(z_{i}+\Delta z)q_{\phi, f}(a(z_{i}) + \Delta a) = Z_e(z_{i+1}+\Delta z)q_{\phi, e}(b(z_{i+1}) + \Delta b)
\end{align}

We now use Taylor expansions for each function, given by:
\begin{align*}
    Z_e(z_{i+1} + \Delta z) &= Z_e(z_{i+1}) + \Delta z Z'_e(z_{i+1}) + ... \\
    Z_f(z_i + \Delta z) &= Z_f(z_i) + \Delta z Z'_f(z_i) + ... \\
    q_{\phi, f}(a(z_i) + \Delta a(z_i)) &= q_{\phi, f}(a(z_i)) + \Delta a(z_i) q'_{\phi, f}(a(z_i)) + ...\\
    q_{\phi, e}(b(z_{i+1}) + \Delta b(z_{i+1})) &= q_{\phi, f}(b(z_{i+1})) + \Delta b(z_{i+1}) q'_{\phi, f}(b(z_{i+1})) + ...\\
    q_{\psi, f}(a(z_i) + \Delta a(z_i)) &= q_{\psi, f}(a(z_i)) + \Delta a(z_i) q'_{\psi, f}(a(z_i)) + ...\\
    q_{\psi, e}(b(z_{i+1}) + \Delta b(z_{i+1})) &= q_{\psi, f}(b(z_{i+1})) + \Delta b(z_{i+1}) q'_{\psi, f}(b(z_{i+1})) + ...
\end{align*}

where $\Delta z$ is chosen such that higher order terms can be ignored. Substituting these expansions into equations (39-40) yields:
\begin{align}
    \big[ Z_f(z_i) + &\Delta z Z'_f(z_i) \big] (q_{\psi, f}(a(z_i)) + \Delta a(z_{i}) q'_{\psi, f}(a(z_i))) \nonumber \\
    & = \left[ Z_e(z_{i+1}) + \Delta z Z'_e(z_{i+1}) \right] (q_{\psi, f}(b(z_{i+1})) + \Delta b(z_{i+1}) q'_{\psi, f}(b(z_{i+1}))) \\
    \nonumber \\
    \big[ Z_f(z_i) + &\Delta z Z'_f(z_i) \big] (q_{\phi, f}(a(z_i)) + \Delta a(z_{i}) q'_{\phi, f}(a(z_i))) \nonumber \\
    & = \left[ Z_e(z_{i+1}) + \Delta z Z'_e(z_{i+1}) \right] (q_{\phi, f}(b(z_{i+1})) + \Delta b(z_{i+1}) q'_{\phi, f}(b(z_{i+1})))
\end{align}

Consider equation (41), expanding out the brackets:
\begin{align*}
    \big[ Z_f(z_i) + &\Delta z Z'_f(z_i) \big] (q_{\psi, f}(a(z_i)) + \Delta a(z_{i}) q'_{\psi, f}(a(z_i))) \nonumber \\
    & = \left[ Z_e(z_{i+1}) + \Delta z Z'_e(z_{i+1}) \right] (q_{\psi, f}(b(z_{i+1})) + \Delta b(z_{i+1}) q'_{\psi, f}(b(z_{i+1}))) \\
\end{align*}
\newpage

Expanding out the brackets:
\begin{align*} 
&Z_f(z_i)q_{\psi, f}(a(z_i)) + \Delta a(z_{i}) Z_f(z_i)q'_{\psi, f}(a(z_i)) + \Delta z Z'_f(z_i)q_{\psi, f}(a(z_i)) +  \\
& \Delta z \Delta a(z_{i}) Z'_f(z_i)q'_{\psi, f}(a(z_i)) = Z_e(z_{i+1})q_{\psi, e}(b(z_{i+1})) + \Delta b(z_{i+1}) Z_e(z_{i+1})q'_{\psi, e}(b(z_{i+1})) + \\
& \Delta z Z'_e(z_{i+1})q_{\psi, e}(b(z_{i+1})) + \Delta z \Delta b(z_{i+1}) Z'_e(z_{i+1})q'_{\psi, e}(b(z_{i+1}))
\intertext{Using (35) to remove the first term from both equations:}
& \Delta a(z_{i}) Z_f(z_i)q'_{\psi, f}(a(z_i)) + \Delta z Z'_f(z_i)q_{\psi, f}(a(z_i)) + \Delta z \Delta a(z_{i}) Z'_f(z_i)q'_{\psi, f}(a(z_i)) = \\
& \Delta b(z_{i+1}) Z_e(z_{i+1})q'_{\psi, e}(b(z_{i+1})) + \Delta z Z'_e(z_{i+1})q_{\psi, e}(b(z_{i+1})) + \Delta z \Delta b(z_{i+1}) Z'_e(z_{i+1})q'_{\psi, e}(b(z_{i+1})) \\
\intertext{Grouping Similar terms:}
&\Delta a(z_{i}) \left( Z_f(z_i) + \Delta z Z'_f(z_i) \right)q'_{\psi, f}(a(z_i)) = \\
&\Delta b(z_{i+1})\left( Z_e(z_{i+1}) + \Delta z Z_e(z_{i+1}) \right)q'_{\phi, e}(b(z_{i+1})) + \Delta z \left( Z'_e(z_{i+1})q_{\phi, e}(b(z_{i+1})) - Z'_f(z_i)q_{\phi,f}(a(z_i)) \right) \\ 
\end{align*}
Rearranging to isolate $\Delta a(z_{i})$
\begin{align}
\Delta a(z_{i}) = &\Delta b(z_{i+1}) \frac{ \left( Z_e(z_{i+1}) + \Delta z Z'_e(z_{i+1}) \right) q'_{\psi,e}(b(z_{i+1}))}{( Z_f(z_i) + \Delta z Z'_f(z_i) ) q'_{\psi, f}(a(z_i))} + \quad \nonumber \\
&\qquad \qquad \qquad \Delta z \frac{( Z'_e(z_{i+1})q_{\psi, e}(b(z_{i+1})) - Z'_f(z_i)q_{\psi, f}(a(z_i)) )}{( Z_f(z_i) + \Delta z Z'_f(z_i) )q'_{\psi, f}(a(z_i))}.
\end{align}

Similarly, if we were to consider equation (42) then we would find the following relationship:
\begin{align}
    \Delta a(z_{i}) = &\Delta b(z_{i+1}) \frac{ \left( Z_e(z_{i+1}) + \Delta z Z'_e(z_{i+1}) \right) q'_{\phi,e}(b(z_{i+1}))}{( Z_f(z_i) + \Delta z Z'_f(z_i) ) q'_{\phi, f}(a(z_i))} \quad + \nonumber \\
    &\qquad \qquad \qquad \Delta z \frac{( Z'_e(z_{i+1})q_{\phi, e}(b(z_{i+1})) - Z'_f(z_i)q_{\phi, f}(a(z_i)) )}{( Z_f(z_i) + \Delta z Z'_f(z_i) )q'_{\phi, f}(a(z_i))}
\end{align}

\newpage

Equating these relationships allows us to eliminate $a(z_i)$, to obtain a relationship between the upper optimal control point and the independent variable $z$:
\begin{align}
    &\Delta b(z_{i+1}) \frac{ \left( Z_e(z_{i+1}) + \Delta z Z'_e(z_{i+1}) \right) q'_{\psi,e}(b(z_{i+1}))}{( Z_f(z_i) + \Delta z Z'_f(z_i) ) q'_{\psi, f}(a(z_i))} + \Delta z \frac{( Z'_e(z_{i+1})q_{\psi, e}(b(z_{i+1})) - Z'_f(z_i)q_{\psi, f}(a(z_i)) )}{( Z_f(z_i) + \Delta z Z'_f(z_i) )q'_{\psi, f}(a(z_i))} \nonumber \\
    &= \Delta b(z_{i+1}) \frac{ \left( Z_e(z_{i+1}) + \Delta z Z'_e(z_{i+1}) \right) q'_{\phi,e}(b(z_{i+1}))}{( Z_f(z_i) + \Delta z Z'_f(z_i) ) q'_{\phi, f}(a(z_i))} + \Delta z \frac{( Z'_e(z_{i+1})q_{\phi, e}(b(z_{i+1})) - Z'_f(z_i)q_{\phi, f}(a(z_i)) )}{( Z_f(z_i) + \Delta z Z'_f(z_i) )q'_{\phi, f}(a(z_i))} \nonumber \\
    \intertext{grouping similar terms and rearranging:}
    &\Delta b(z_{i+1}) \left[ \frac{ \left( Z_e(z_{i+1}) + \Delta z Z'_e(z_{i+1}) \right) q'_{\phi,e}(b(z_{i+1}))}{( Z_f(z_i) + \Delta z Z'_f(z_i) ) q'_{\phi, f}(a(z_i))} - \frac{ \left( Z_e(z_{i+1}) + \Delta z Z'_e(z_{i+1}) \right) q'_{\psi,e}(b(z_{i+1}))}{( Z_f(z_i) + \Delta z Z'_f(z_i) ) q'_{\psi, f}(a(z_i))} \right] \nonumber \\
    &= \Delta z \left[ \frac{( Z'_e(z_{i+1})q_{\psi, e}(b(z_{i+1})) - Z'_f(z_i)q_{\psi, f}(a(z_i)) )}{( Z_f(z_i) + \Delta z Z'_f(z_i) )q'_{\psi, f}(a(z_i))} - \frac{( Z'_e(z_{i+1})q_{\phi, e}(b(z_{i+1})) - Z'_f(z_i)q_{\phi, f}(a(z_i)) )}{( Z_f(z_i) + \Delta z Z'_f(z_i) )q'_{\phi, f}(a(z_i))} \right] \nonumber \\
    \intertext{Removing the common factor of $\frac{1}{Z_f(z_i)+\Delta z Z'_f(z_i)}$, and taking $(Z_e(z_{i+1})+\Delta z Z'_e(z_{i+1}))$ outside the LHS}
    \implies &\Delta b(z_{i+1}) \left( Z_e(z_{i+1}) + \Delta z Z'_e(z_{i+1}) \right) \left[ \frac{  q'_{\phi,e}(b(z_{i+1}))}{ q'_{\phi, f}(a(z_i))} - \frac{  q'_{\psi,e}(b(z_{i+1}))}{ q'_{\psi, f}(a(z_i))} \right] \nonumber \\
    &= \Delta z \left[ \frac{( Z'_e(z_{i+1})q_{\psi, e}(b(z_{i+1})) - Z'_f(z_i)q_{\psi, f}(a(z_i)) )}{q'_{\psi, f}(a(z_i))} - \frac{( Z'_e(z_{i+1})q_{\phi, e}(b(z_{i+1})) - Z'_f(z_i)q_{\phi, f}(a(z_i)) )}{q'_{\phi, f}(a(z_i))} \right] \nonumber \\
    \intertext{Obtaining a common denominator on the RHS}
    &= \Delta z \Bigg[ \frac{ Z'_e(z_{i+1})q'_{\phi, f}(a(z_i))q_{\psi, e}(b(z_{i+1})) - Z'_e(z_{i+1})q'_{\psi, f}(a(z_i))q_{\phi, e}(b(z_{i+1})) }{q'_{\psi, f}(a(z_i))q'_{\phi, f}(a(z_i))} \nonumber \\
    &\quad \quad \quad \quad + \frac{- Z'_f(z_i)q'_{\phi, f}(a(z_i))q_{\psi, f}(a(z_i))  + Z'_f(z_i)q_{\phi, f}(a(z_i))q'_{\psi, f}(a(z_i))}{q'_{\psi, f}(a(z_i))q'_{\phi, f}(a(z_i))} \Bigg] \nonumber \\
    \intertext{Grouping terms by $Z$ functions:}
    &= \Delta z \Bigg[ \frac{Z'_e(z_{i+1}) \left( q_{\psi,e}(b(z_{i+1}))q'_{\phi,f}(a(z_i)) - q'_{\psi,f}(a(z_i))q_{\phi,e}(b(z_{i+1})) \right)}{q'_{\psi, f}(a(z_i))q'_{\phi, f}(a(z_i))} \nonumber \\
    &\quad \quad \quad \quad + \frac{Z'_f(z_i) \left(q_{\phi, f}(a(z_i))q'_{\psi, f}(a(z_i)) - q'_{\phi, f}(a(z_i))q_{\psi, f}(a(z_i)) \right)}{q'_{\psi, f}(a(z_i))q'_{\phi, f}(a(z_i))} \Bigg] \nonumber \\
    \intertext{Simplifying gives:}
    &= \Delta z \Bigg[ Z'_e(z_{i+1})\left(\frac{q_{\psi,e}(b(z_{i+1}))}{q'_{\psi,f}(a(z_{i}))} - \frac{q_{\phi,e}(b(z_{i+1}))}{q'_{\phi,f}(a(z_{i}))}\right) + Z'_f(z_i) \left( \frac{q_{\phi, e}(a(z_i)}{q'_{\phi, e}(a(z_i)} - \frac{q_{\psi, e}(a(z_i)}{q'_{\psi, e}(a(z_i)} \right) \Bigg] \nonumber \\
    \intertext{Rearranging gives the following result for $\frac{\Delta b(z_{i+1})}{\Delta z}$}
    \implies & \frac{\Delta b(z_{i+1})}{\Delta z} =  \frac{\Bigg[ Z'_e(z_{i+1})\left(\frac{q_{\psi,e}(b(z_{i+1}))}{q'_{\psi,f}(a(z_{i}))} - \frac{q_{\phi,e}(b(z_{i+1}))}{q'_{\phi,f}(a(z_{i}))}\right) + Z'_f(z_i) \left( \frac{q_{\phi, e}(a(z_i)}{q'_{\phi, e}(a(z_i)} - \frac{q_{\psi, e}(a(z_i)}{q'_{\psi, e}(a(z_i)} \right) \Bigg]}{\left( Z_e(z_{i+1}) + \Delta z Z'_e(z_{i+1}) \right)\left[ \frac{  q'_{\phi,e}(b(z_{i+1}))}{ q'_{\phi, f}(a(z_i))} - \frac{q'_{\psi,e}(b(z_{i+1}))}{ q'_{\psi, f}(a(z_i))} \right]}.
\end{align}

We look to simplify this expression, starting with the $Z'_e(z_{i+1})$ term:
\begin{align*}
    Z'_e(z_{i+1})\frac{\left[\frac{q_{\psi,e}(b(z_{i+1}))}{q'_{\psi,f}(a(z_{i}))} - \frac{q_{\phi,e}(b(z_{i+1}))}{q'_{\phi,f}(a(z_{i}))}\right]}{\left[ \frac{  q'_{\phi,e}(b(z_{i+1}))}{ q'_{\phi, f}(a(z_i))} - \frac{q'_{\psi,e}(b(z_{i+1}))}{ q'_{\psi, f}(a(z_i))} \right]} &= Z'_e(z_{i+1}) \frac{\frac{q_{\psi, e}(b(z_{i+1}))q'_{\phi,f}(a(z_i)) - q'_{\psi,f}(a(z_i))q_{\phi,e}(b(z_{i+1}))}{q'_{\psi,f}(a(z_i))q'_{\phi,f}(a(z_i))}}{\frac{q'_{\psi,f}(a(z_i))q'_{\phi,e}(b(z_{i+1})) - q'_{\psi,e}(b(z_{i+1}))q'_{\phi,f}(a(z_i))}{q'_{\psi,f}(a(z_i))q'_{\phi,f}(a(z_i))}} \\
    &= Z'_e(z_{i})\frac{q_{\psi, e}(b(z_{i+1}))q'_{\phi,f}(a(z_i)) - q'_{\psi,f}(a(z_i))q_{\phi,e}(b(z_{i+1}))}{q'_{\psi,f}(a(z_i))q'_{\phi,e}(b(z_{i+1})) - q'_{\psi,e}(b(z_{i+1}))q'_{\phi,f}(a(z_i))}
\end{align*}
then considering the $Z'_f(z_i)$ term:
\begin{align*}
    Z'_f(z_i) \frac{\left[ \frac{q_{\phi, e}(a(z_i))}{q'_{\phi, e}(a(z_i))} - \frac{q_{\psi, e}(a(z_i))}{q'_{\psi, e}(a(z_i))} \right]}{\left[ \frac{  q'_{\phi,e}(b(z_{i+1}))}{q'_{\phi, f}(a(z_i))} - \frac{q'_{\psi,e}(b(z_{i+1}))}{ q'_{\psi, f}(a(z_i))} \right]} &= Z'_f(z_i) \frac{\frac{q'_{\psi,f}(a(z_i))q_{\phi,f}(a(z_i)) - q_{\psi,f}(a(z_i))q'_{\phi,f}(a(z_i))}{q'_{\psi,f}(a(z_i))q'_{\phi,f}(a(z_i))}}{\frac{q'_{\psi,f}(a(z_i))q'_{\phi,e}(b(z_{i+1})) - q'_{\psi,f}(a(z_i))q'_{\phi,f}(a(z_i))}{q'_{\psi,f}(a(z_i))q'_{\phi,f}(a(z_i))}} \\
    &= Z'_f(z_{i})\frac{q'_{\psi,f}(a(z_i))q_{\phi,f}(a(z_i)) - q_{\psi,f}(a(z_i))q'_{\phi,f}(a(z_i))}{q'_{\psi,f}(a(z_i))q'_{\phi,e}(b(z_{i+1})) - q'_{\psi,e}(b(z_{i+1}))q'_{\phi,f}(a(z_i))}
\end{align*}
Substituting these results into (45) we find:
\begin{align*}
        \frac{\Delta b(z_{i+1})}{\Delta z} = \frac{1}{Z_e(z_{i+1}) + \Delta z Z'_e(z_{i+1})}&\Bigg( Z'_e(z_{i+1})\left[\frac{ q_{\psi,e}(b(z_{i+1}))q'_{\phi,f}(a(z_i)) - q'_{\psi,f}(a(z_i))q_{\phi,e}(b(z_{i+1})) }{q'_{\psi,f}(a(z_i))q'_{\phi,e}(b(z_{i+1})) - q'_{\psi,e}(b(z_{i+1}))q'_{\phi,f}(a(z_i))}\right] \\
        & \quad + Z'_f(z_i)\left[\frac{ q'_{\psi,f}(a(z_i))q_{\phi,f}(a(z_i)) - q_{\psi,f}(a(z_i))q'_{\phi,f}(a(z_i)) }{q'_{\psi,f}(a(z_i))q'_{\phi,e}(b(z_{i+1})) - q'_{\psi,e}(b(z_{i+1}))q'_{\phi,f}(a(z_i))}\right] \Bigg).
\end{align*}

\newpage
So we have now obtained an expression for $\Delta b(z_{i+1})$, and we can substitute this back into equation (43):
\begin{align*}
    \frac{\Delta a(z_i)}{\Delta z} &= \frac{1}{Z_f(z_i) + \Delta z Z'_f(z_i)}\Bigg( Z'_e(z_{i+1})\left[\frac{ q_{\psi,e}(b(z_{i+1}))q'_{\phi,f}(a(z_i)) - q'_{\psi,f}(a(z_i))q_{\phi,e}(b(z_{i+1})) }{q'_{\psi,f}(a(z_i))q'_{\phi,e}(b(z_{i+1})) - q'_{\psi,e}(b(z_{i+1}))q'_{\phi,f}(a(z_i))}\right] \\
    & \quad + Z'_f(z_i)\left[\frac{ q'_{\psi,f}(a(z_i))q_{\phi,f}(a(z_i)) - q_{\psi,f}(a(z_i))q'_{\phi,f}(a(z_i)) }{q'_{\psi,f}(a(z_i))q'_{\phi,e}(b(z_{i+1})) - q'_{\psi,e}(b(z_{i+1}))q'_{\phi,f}(a(z_i))}\right] \Bigg)\frac{q'_{\phi,e}(b(z_{i+1}))}{q'_{\phi,f}(a(z_i))} \\
    & \quad + \frac{( Z'_e(z_{i+1})q_{\phi, e}(b(z_{i+1})) - Z'_f(z_i)q_{\phi, f}(a(z_i)) )}{( Z_f(z_i) + \Delta z Z'_f(z_i) )q'_{\phi, f}(a(z_i))} \\
    \intertext{Factorising the denominators}
    &= \frac{1}{Z_f(z_i) + \Delta z Z'_f(z_i)}\frac{1}{q'_{\phi,f}(a(z_i))}\Bigg[ \Bigg( Z'_e(z_{i+1})\left[\frac{ q_{\psi,e}(b(z_{i+1}))q'_{\phi,f}(a(z_i)) - q'_{\psi,f}(a(z_i))q_{\phi,e}(b(z_{i+1})) }{q'_{\psi,f}(a(z_i))q'_{\phi,e}(b(z_{i+1})) - q'_{\psi,e}(b(z_{i+1}))q'_{\phi,f}(a(z_i))}\right] \\
    & \quad + Z'_f(z_i)\left[\frac{ q'_{\psi,f}(a(z_i))q_{\phi,f}(a(z_i)) - q_{\psi,f}(a(z_i))q'_{\phi,f}(a(z_i)) }{q'_{\psi,f}(a(z_i))q'_{\phi,e}(b(z_{i+1})) - q'_{\psi,e}(b(z_{i+1}))q'_{\phi,f}(a(z_i))}\right] \Bigg)q'_{\phi,e}(b(z_{i+1})) \\
    & \quad + ( Z'_e(z_{i+1})q_{\phi, e}(b(z_{i+1})) - Z'_f(z_i)q_{\phi, f}(a(z_i)) ) \Bigg]
    \intertext{Gathering similar terms}
    &= \frac{1}{Z_f(z_i) + \Delta z Z'_f(z_i)} \frac{1}{q'_{\phi,f}(a(z_i))} \\
    & \quad \quad \Bigg( Z'_e(z_{i+1})\left[\frac{ q_{\psi,e}(b(z_{i+1}))q'_{\phi,f}(a(z_i)) - q'_{\psi,f}(a(z_i))q_{\phi,e}(b(z_{i+1})) }{q'_{\psi,f}(a(z_i))q'_{\phi,e}(b(z_{i+1})) - q'_{\psi,e}(b(z_{i+1}))q'_{\phi,f}(a(z_i))}q'_{\phi,e}(b(z_{i+1})) + q_{\phi,e}(b(z_{i+1}))\right] \\
    & \quad \quad + Z'_f(z_i)\left[\frac{ q'_{\psi,f}(a(z_i))q_{\phi,f}(a(z_i)) - q_{\psi,f}(a(z_i))q'_{\phi,f}(a(z_i)) }{q'_{\psi,f}(a(z_i))q'_{\phi,e}(b(z_{i+1})) - q'_{\psi,e}(b(z_{i+1}))q'_{\phi,f}(a(z_i))}q'_{\phi,e}(b(z_{i+1})) - q_{\phi,f}(a(z_i))  \right] \Bigg).
\end{align*}

We start simplifying by considering the $Z'_e(z_{i+1})$ term:
\begin{align*}
    &\left[\frac{ q_{\psi,e}(b(z_{i+1}))q'_{\phi,f}(a(z_i)) - q'_{\psi,f}(a(z_i))q_{\phi,e}(b(z_{i+1})) }{q'_{\psi,f}(a(z_i))q'_{\phi,e}(b(z_{i+1})) - q'_{\psi,e}(b(z_{i+1}))q'_{\phi,f}(a(z_i))}q'_{\phi,e}(b(z_{i+1})) + q_{\phi,e}(b(z_{i+1}))\right] \\
    \intertext{First equalising the denominators}
    &= \Bigg[\frac{ q_{\psi,e}(b(z_{i+1}))q'_{\phi,f}(a(z_i))q'_{\phi,e}(b(z_{i+1})) - q'_{\psi,f}(a(z_i))q_{\phi,e}(b(z_{i+1}))q'_{\phi,e}(b(z_{i+1})) }{q'_{\psi,f}(a(z_i))q'_{\phi,e}(b(z_{i+1})) - q'_{\psi,e}(b(z_{i+1}))q'_{\phi,f}(a(z_i))} \\
    & \quad \quad + \frac{q'_{\psi,f}(a(z_i))q_{\phi,e}(b(z_{i+1}))q'_{\phi,e}(b(z_{i+1})) - q'_{\psi,e}(b(z_{i+1}))q'_{\phi,f}(a(z_i))q_{\phi,e}(b(z_{i+1}))}{q'_{\psi,f}(a(z_i))q'_{\phi,e}(b(z_{i+1})) - q'_{\psi,e}(b(z_{i+1}))q'_{\phi,f}(a(z_i))}\Bigg] \\ 
    \intertext{Then cancelling equivalent terms}
    &= \Bigg[\frac{ q_{\psi,e}(b(z_{i+1}))q'_{\phi,f}(a(z_i))q'_{\phi,e}(b(z_{i+1})) - q'_{\psi,e}(b(z_{i+1}))q'_{\phi,f}(a(z_i))q_{\phi,e}(b(z_{i+1}))}{q'_{\psi,f}(a(z_i))q'_{\phi,e}(b(z_{i+1})) - q'_{\psi,e}(b(z_{i+1}))q'_{\phi,f}(a(z_i))} \Bigg] \\
    \intertext{Finally factorising gives}
    &= q'_{\phi,f}(a(z_i))\Bigg[\frac{ q_{\psi,e}(b(z_{i+1}))q'_{\phi,e}(b(z_{i+1})) - q'_{\psi,e}(b(z_{i+1}))q_{\phi,e}(b(z_{i+1}))}{q'_{\psi,f}(a(z_i))q'_{\phi,e}(b(z_{i+1})) - q'_{\psi,e}(b(z_{i+1}))q'_{\phi,f}(a(z_i))} \Bigg]
\end{align*}
Then the $Z'_f(z_i)$ term:
\begin{align*}
    &\left[\frac{ q'_{\psi,f}(a(z_i))q_{\phi,f}(a(z_i)) - q_{\psi,f}(a(z_i))q'_{\phi,f}(a(z_i)) }{q'_{\psi,f}(a(z_i))q'_{\phi,e}(b(z_{i+1})) - q'_{\psi,e}(b(z_{i+1}))q'_{\phi,f}(a(z_i))}q'_{\phi,e}(b(z_{i+1})) - q_{\phi,f}(a(z_i))  \right] \\
    \intertext{First equalising the denominators}
    &= \Bigg[\frac{ q'_{\psi,f}(a(z_i))q_{\phi,f}(a(z_i))q'_{\phi,e}(b(z_{i+1})) - q_{\psi,f}(a(z_i))q'_{\phi,f}(a(z_i))q'_{\phi,e}(b(z_{i+1})) }{q'_{\psi,f}(a(z_i))q'_{\phi,e}(b(z_{i+1})) - q'_{\psi,e}(b(z_{i+1}))q'_{\phi,f}(a(z_i))} \\
    &\quad \quad + \frac{-q'_{\psi,f}(a(z_i))q_{\phi,f}(a(z_i))q'_{\phi,e}(b(z_{i+1})) + q'_{\psi,e}(b(z_{i+1}))q_{\phi,f}(a(z_i))q'_{\phi,f}(a(z_i))}{q'_{\psi,f}(a(z_i))q'_{\phi,e}(b(z_{i+1})) - q'_{\psi,e}(b(z_{i+1}))q'_{\phi,f}(a(z_i))}  \Bigg] \\
    \intertext{Then cancelling equivalent terms}
    &= \Bigg[\frac{ q'_{\psi,e}(b(z_{i+1}))q_{\phi,f}(a(z_i))q'_{\phi,f}(a(z_i)) - q_{\psi,f}(a(z_i))q'_{\phi,f}(a(z_i))q'_{\phi,e}(b(z_{i+1})) }{q'_{\psi,f}(a(z_i))q'_{\phi,e}(b(z_{i+1})) - q'_{\psi,e}(b(z_{i+1}))q'_{\phi,f}(a(z_i))} \Bigg] \\
    \intertext{Finally factorising gives}
    &= q'_{\phi,f}(a(z_i)) \Bigg[\frac{ q'_{\psi,e}(b(z_{i+1}))q_{\phi,f}(a(z_i)) - q_{\psi,f}(a(z_i))q'_{\phi,e}(b(z_{i+1})) }{q'_{\psi,f}(a(z_i))q'_{\phi,e}(b(z_{i+1})) - q'_{\psi,e}(b(z_{i+1}))q'_{\phi,f}(a(z_i))} \Bigg]
\end{align*}
Substituting these results into our expression for $\frac{\Delta a(z_i)}{\Delta z}$:
\begin{align*}
    \frac{\Delta a(z_i)}{\Delta z} &= \frac{1}{Z_f(z_i) + \Delta z Z'_f(z_i)} \frac{1}{q'_{\phi,f}(a(z_i))} \\
    & \quad \quad \Bigg( Z'_e(z_{i+1})\left[ q'_{\phi,f}(a(z_i))\Bigg(\frac{ q_{\psi,e}(b(z_{i+1}))q'_{\phi,e}(b(z_{i+1})) - q'_{\psi,e}(b(z_{i+1}))q_{\phi,e}(b(z_{i+1}))}{q'_{\psi,f}(a(z_i))q'_{\phi,e}(b(z_{i+1})) - q'_{\psi,e}(b(z_{i+1}))q'_{\phi,f}(a(z_i))} \Bigg) \right] \\
    & \quad \quad + Z'_f(z_i)\Bigg[ q'_{\phi,f}(a(z_i)) \Bigg(\frac{ q'_{\psi,e}(b(z_{i+1}))q_{\phi,f}(a(z_i)) - q_{\psi,f}(a(z_i))q'_{\phi,e}(b(z_{i+1})) }{q'_{\psi,f}(a(z_i))q'_{\phi,e}(b(z_{i+1})) - q'_{\psi,e}(b(z_{i+1}))q'_{\phi,f}(a(z_i))} \Bigg) \Bigg]
    \Bigg).
\end{align*}

After cancelling, we find the following expression for the derivative of our lower optimal control point with the independent variable $z$, alongside our previous expression for the upper optimal control point:
\begin{align*}
    \frac{\Delta a(z_i)}{\Delta z} = \frac{1}{Z_f(z_i) + \Delta z Z'_f(z_i)}\Bigg( &Z'_e(z_{i+1})\left[ \Bigg(\frac{ q_{\psi,e}(b(z_{i+1}))q'_{\phi,e}(b(z_{i+1})) - q'_{\psi,e}(b(z_{i+1}))q_{\phi,e}(b(z_{i+1}))}{q'_{\psi,f}(a(z_i))q'_{\phi,e}(b(z_{i+1})) - q'_{\psi,e}(b(z_{i+1}))q'_{\phi,f}(a(z_i))} \Bigg) \right] \\
    + &Z'_f(z_i)\Bigg[  \Bigg(\frac{ q'_{\psi,e}(b(z_{i+1}))q_{\phi,f}(a(z_i)) - q_{\psi,f}(a(z_i))q'_{\phi,e}(b(z_{i+1})) }{q'_{\psi,f}(a(z_i))q'_{\phi,e}(b(z_{i+1})) - q'_{\psi,e}(b(z_{i+1}))q'_{\phi,f}(a(z_i))} \Bigg) \Bigg] \Bigg)
\end{align*}

\begin{align*}
        \frac{\Delta b(z_{i+1})}{\Delta z} = \frac{1}{Z_e(z_{i+1}) + \Delta z Z'_e(z_{i+1})}\Bigg( &Z'_e(z_{i+1})\left[\frac{ q_{\psi,e}(b(z_{i+1}))q'_{\phi,f}(a(z_i)) - q'_{\psi,f}(a(z_i))q_{\phi,e}(b(z_{i+1})) }{q'_{\psi,f}(a(z_i))q'_{\phi,e}(b(z_{i+1})) - q'_{\psi,e}(b(z_{i+1}))q'_{\phi,f}(a(z_i))}\right] \\
        + &Z'_f(z_i)\left[\frac{ q'_{\psi,f}(a(z_i))q_{\phi,f}(a(z_i)) - q_{\psi,f}(a(z_i))q'_{\phi,f}(a(z_i)) }{q'_{\psi,f}(a(z_i))q'_{\phi,e}(b(z_{i+1})) - q'_{\psi,e}(b(z_{i+1}))q'_{\phi,f}(a(z_i))}\right] \Bigg).
\end{align*}

We can now take the limit as $\Delta z \rightarrow 0$ to derive expressions for the derivatives of the optimal control points with regards to the independent variable:

\begin{align*}
\lim_{\Delta z \rightarrow 0} \frac{\Delta a(z_i)}{\Delta z} = \frac{d(a(z))}{dz} = \frac{1}{Z_f(z)}\Bigg( &Z'_e(z)\left[ \Bigg(\frac{ q_{\psi,e}(b(z))q'_{\phi,e}(b(z)) - q'_{\psi,e}(b(z))q_{\phi,e}(b(z))}{q'_{\psi,f}(a(z))q'_{\phi,e}(b(z)) - q'_{\psi,e}(b(z))q'_{\phi,f}(a(z))} \Bigg) \right] \\
+ &Z'_f(z)\Bigg[  \Bigg(\frac{ q'_{\psi,e}(b(z))q_{\phi,f}(a(z)) - q_{\psi,f}(a(z))q'_{\phi,e}(b(z)) }{q'_{\psi,f}(a(z))q'_{\phi,e}(b(z)) - q'_{\psi,e}(b(z))q'_{\phi,f}(a(z))} \Bigg) \Bigg] \Bigg)
\end{align*}

\begin{align*}
\lim_{\Delta z \rightarrow 0} \frac{\Delta b(z_{i+1})}{\Delta z} = \frac{d(b(z))}{dz} = \frac{1}{Z_e(z)}\Bigg( &Z'_e(z)\left[\frac{ q_{\psi,e}(b(z))q'_{\phi,f}(a(z)) - q'_{\psi,f}(a(z))q_{\phi,e}(b(z)) }{q'_{\psi,f}(a(z))q'_{\phi,e}(b(z)) - q'_{\psi,e}(b(z))q'_{\phi,f}(a(z))}\right] \\
+ &Z'_f(z)\left[\frac{ q'_{\psi,f}(a(z))q_{\phi,f}(a(z)) - q_{\psi,f}(a(z))q'_{\phi,f}(a(z)) }{q'_{\psi,f}(a(z))q'_{\phi,e}(b(z)) - q'_{\psi,e}(b(z))q'_{\phi,f}(a(z))}\right] \Bigg).
\end{align*}

And finally use the limit of equations (35-36) taken as $\Delta z \rightarrow 0$ to find the final results:

\begin{align*}
    \frac{da}{dz} &= \left[ \frac{Z'_e(z)}{Z_e(z)} - \frac{Z'_f(z)}{Z_f(z)} \right] \left( \frac{q_{\psi,f}(a)q'_{\phi,e}(b) - q'_{\psi,e}(b)q_{\phi,f}(a)}{q'_{\psi,f}(a)q'_{\phi,e}(b) - q'_{\psi,e}(b)q'_{\phi,f}(a)} \right) \\
    &= \frac{d}{dz}\left( \frac{Z_e(z)}{Z_f(z)} \right)\left( \frac{q_{\psi,e}(b)q'_{\phi,e}(b) - q'_{\psi,e}(b)q_{\phi,e}(b)}{q'_{\psi,f}(a)q'_{\phi,e}(b) - q'_{\psi,e}(b)q'_{\phi,f}(a)} \right) \numberthis
\end{align*}

\begin{align*}
    \frac{db}{dz} &= \left[ \frac{Z'_e(z)}{Z_e(z)} - \frac{Z'_f(z)}{Z_f(z)} \right] \left( \frac{q_{\psi,e}(b)q'_{\phi,f}(a) - q'_{\psi,f}(a)q_{\phi,e}(b)}{q'_{\psi,f}(a)q'_{\phi,e}(b) - q'_{\psi,e}(b)q'_{\phi,f}(a)} \right) \\
    &= \frac{d}{dz}\left( \frac{Z_f(z)}{Z_e(z)} \right)\left( \frac{q_{\psi,f}(a)q'_{\phi,f}(a) - q'_{\psi,f}(a)q_{\phi,f}(a)}{q'_{\psi,f}(a)q'_{\phi,e}(b) - q'_{\psi,e}(b)q'_{\phi,f}(a)} \right) \numberthis
\end{align*}

\end{document}